\documentclass[12pt]{article}
\usepackage[cp866]{inputenc}
\usepackage[english]{babel}
\usepackage{graphicx}
\usepackage{subfig}
\usepackage{extsizes}
\usepackage{epsfig}
\usepackage[numbers,sort&compress]{natbib}
\usepackage{citehack}
\usepackage{float}
\usepackage{amsmath}
\usepackage{amssymb}
\usepackage{multicol}
\usepackage{amscd}
\usepackage{verbatim}
\usepackage{booktabs}
\usepackage{enumitem}
\newcommand{\eps}{\varepsilon}
\newcommand{\vte}{\vartheta}
\newcommand{\vfi}{\varphi}

\newcommand{\R}{\mathbb {R}}

\newcommand{\Sm}{\mathbb {S}}
\newcommand{\p}{\partial}

\newtheorem{theorem}{Theorem}
\newtheorem{proposition}{Proposition}
\newtheorem{lemma}{Lemma}[section]

\newtheorem{dfn}{Definition}

\makeatletter \@addtoreset{equation}{section} \makeatother
\allowdisplaybreaks[1]

\begin{document}
      \title{On phase at a resonance in slow-fast Hamiltonian systems }
      \date{}
      \author{Yuyang Gao,  Anatoly Neishtadt, Alexey Okunev} 
      \maketitle
    \begin{abstract}
    We consider a slow-fast Hamiltonian system with one fast angular variable (a fast phase)  whose frequency vanishes on some surface  in the space of slow variables (a resonant surface). Systems of such form appear in the study of dynamics of charged particles in inhomogeneous  magnetic field under influence of a high-frequency electrostatic waves. Trajectories of the averaged over the fast phase system cross the resonant surface. The fast phase makes $\sim \frac 1\eps$ turns before arrival to the resonant  surface ($\eps$ is a small parameter of the problem). An asymptotic formula for the value of the phase at the arrival to the resonance was derived earlier in the context of study of charged particle dynamics on the basis of heuristic considerations  without any estimates of its accuracy. We provide a rigorous derivation of this formula and prove that its accuracy is $O(\sqrt \eps)$ (up to a logarithmic correction). Numerics indicate that this estimate for the accuracy is optimal.
    \end{abstract}

      \section{Introduction}

      Slow-fast dynamical systems are systems that contain a small parameter,  say, $\eps$, and variables of two types: slow, changing with a rate $\sim\eps$, and fast, changing with a rate $\sim 1$. Usually, the dynamics of slow variables is of the principal interest. However, there are systems such that the character of the slow dynamics is determined by the dynamics of fast variables. In particular, this is the case when fast variables are angular variables (rotating phases).
      Frequencies of rotation of fast phases depend on slow variables. Evolution  of  the values of slow variables leads to passages through resonances between these frequencies.  Passage through a single resonance leads to a jump in the values of slow variables, and the value of this jump depends on a certain linear combination of the values of the phases
      (\cite {Chir1959, Kevor1974}, see also \cite {Nei1997} and references therein).
      The frequency of this combination of the phases vanishes at the considered resonance. This combination of the phases can be chosen as a new phase for the study of the dynamics near the considered resonance.
      Thus, the underlying phenomenon is a jump of slow variables when the frequency of one of the phases vanishes on some surface in the space of slow variables (a resonant surface). In this paper we consider a model problem when in a slow-fast Hamiltonian system there is just one fast phase.  The frequency of this phase vanishes at some  surface in the space  of slow variables. This is the resonance in the considered problem. An asymptotic formula that expresses the value of the phase at the resonance via initial conditions taken far from the resonance was suggested in the Appendix to \cite{mapping} on the basis of a heuristic consideration without an estimate of the accuracy of this formula. In the current paper we present a rigorous derivation of this formula with an estimate of its accuracy. Phase at a resonance in systems with dissipation is discussed in \cite{Bro_Hab}, \cite{Kis_Tar}.

      Systems of the considered form appear, in particular, in the study of dynamics of charged particles in inhomogeneous magnetic field under influence of a high-frequency electrostatic wave (\cite{mapping} and references therein). The fast phase in these problems is the phase of the wave. The resonance occurs when the projection of the velocity of the particle onto the direction of wave propagation is equal to the phase velocity of the wave. Passage through this resonance leads to a jump in the particle  energy. The value of this jump depends on the phase of the wave at the resonance. The formula for this phase together with the formula for the jump of energy are used in \cite{mapping} to derive a map describing the changes of the energy of the particle and the phase of the wave for consecutive passages through the resonance. The iterations of this map are used in \cite{mapping} to model the long-term dynamics of particles with multiple passages through the resonance. Such an approach allows for a drastic reduction of the modelling time in comparison with direct tracking of particles trajectories by solving the differential equations of the motion numerically.

      \medskip

      The structure of this paper is as follows. In Section \ref{s_outline} we state  the considered problem. The required assumptions are formulated in Section \ref{s_assumptions}. The main result of the paper, the formula for the phase at a resonance (Theorem \ref{theo_1}) is stated in Section \ref {s_main_formula}.  The rest of the paper contains the proof of this formula. In the Appendix we discuss a probabilistic aspect of the problem.

      \section {Outline of the problem}
      \label{s_outline}
      We consider a slow-fast Hamiltonian system with two pairs of canonically conjugate variables
      $ (I,\vfi)$,  $(\eps^{-1}x,y)$  and a Hamilton's function
         \begin{equation}
         \label{H_exact}
        H(I,\vfi,y,x,\eps )=H_0(I,y,x )+\eps H_1(I,\vfi,y,x,\eps ).
         \end{equation}
         Here $\eps>0$ is a small parameter of the problem, $I\in  \R$,  $\vfi\in {\Sm}^1$, $(y,x)\in \R^{2n}, \ n\ge1 $, and $H_1$ is $2\pi$-periodic in $\vfi$.
         The variables $I$ and $\vfi$ are called the \emph{action variable} and the \emph{angle variable} (respectively) of this system.
         The angle variable is also called the \emph{phase}.
      The corresponding differential equations are
      \begin{equation}
      \begin{aligned}
      \label{e_system}
      \dot I&=-\eps \frac{\p H_1}{\p \vfi}, \quad  \dot \vfi= \frac{\p H_0}{\p I}+\eps \frac{\p H_1}{\p I},\\
       \dot y&= -\eps\frac{\p H_0}{\p x}-\eps^2 \frac{\p H_1}{\p x},\quad  \dot x= \eps\frac{\p H_0}{\p y}+\eps^2 \frac{\p H_1}{\p y}\,.
       \end{aligned}
\end{equation}
Thus, $I, y, x$ are slow variables, and $\vfi$ is a fast variable (the fast phase). The function
\begin{equation}
\omega_0(I, y,x)=\frac{\p H_0(I,y,x)}{\p I}
\end {equation}
is the frequency of the fast phase. We assume that this frequency vanishes on some surface ${\cal R}$ in the space of  slow variables. This surface will be called the {\it  resonant surface}.

For an approximate description of the behaviour of the slow variables one can use the {\it averaged system} \cite{bogol}. In the considered case the averaged system has the following simple form
\begin{eqnarray}
      \label{a_system}
      \dot {\bar I}=0, \quad
       \dot {\bar y}= -\eps\frac{\p H_0(\bar I, \bar y, \bar x)}{\p x},\quad  \dot {\bar x}= \eps\frac{\p H_0(\bar I, \bar y, \bar x)}{\p y}\, .
\end{eqnarray}
Thus, for the variables $\bar y, \bar x$
considered as canonically conjugate variables
 we have a Hamiltonian system with the Hamiltonian $ \eps H_0(\bar I, \bar y, \bar x)$ that depends on the parameter $\bar I={\rm const}$. An approximation provided by these formulas is called {\it an adiabatic approximation}, the variable $I$ is {\it an adiabatic invariant}.

Consider a solution $(I(t), \vfi(t), y(t), x(t))$  of the system (\ref{e_system}) with the initial data $(I_0,\vfi_0,  y_0, x_0)$ at $t=t_0$.
Consider the solution $(\bar I,  \bar y(\tau), \bar x(\tau))$ of the averaged system (\ref {a_system}) with the same initial condition for slow variables: $\bar I=I_0$,
  $\bar y(\tau_{0})=y_0$, $\bar x(\tau_0)=x_0$. Here we have introduced the {\it slow time }  $\tau=\eps t$, and $\tau_0=\eps t_0$.
Assume that the considered solution of the averaged system arrives to the resonant surface  for the first time at a moment of slow time $\tau_{*}$. Then, if some mild conditions are satisfied, the considered solution of the exact system also arrives to the resonant surface, suppose that this happens for $t = t_e$. The value of $\eps t_e=\tau_e$ is close to $\tau _{*}$. Let $\vfi_e$ denote the value of the fast phase at the resonance, $\vfi_e=\vfi(t_e)$.
The goal of this paper is to provide  a  rigorous derivation of an asymptotic formula for the value $\vfi_e$ and obtain an estimate of its accuracy.

 \section {Formulation of assumptions}
 \label{s_assumptions}
 System (\ref{e_system}) is considered for $(I, y,x, \vfi )\in D\times \Sm^1$,
 $0 \le \eps\le \eps_1$, where $D$ is a compact domain, $\eps_1={\rm const}$.
Solutions of system (\ref{e_system}) are considered for initial conditions in the domain $D_0\times \Sm^1$, where $D_0$ is a compact domain
 with a piece-wise smooth boundary,
 $D_0\subset D$.   We assume that the following conditions are satisfied.

 \begin{enumerate}[label=\Alph*)]
   \item Functions $H_0, H_1$  are three times continuously  differentiable in $I,\vfi, y,x$. Function $H_1$ is one time continuously differentiable in $\eps$.

   \item  Function $H_0$ is non-degenerate on  the resonant surface ${\cal R}$:
    \begin {equation}
    \frac{\p^2 H_0(I,y,x)}{\p I^2}\equiv \frac{\p\omega_0(I,y,x)}{\p I}\ne 0 \ {\rm at} \ (I,y,x)\in {\cal R}.
    \end{equation}

    \item  Solutions of the averaged system cross the resonant surface ${\cal R}$ transversally:
     \begin {equation}
     \label{transvers}
     \frac{\p\omega_0(I,y,x)}{\p x} \frac{\p H_0(I,y,x)}{\p y}-\frac{\p\omega_0(I,y,x)}{\p y} \frac{\p H_0(I,y,x)}{\p x}\ne 0
     \ {\rm at} \ (I,y,x)\in {\cal R}.
     \end{equation}

     \item  Solutions of the averaged system with  initial conditions in $D_0$ are well defined, stay at a positive distance from the boundary of $D$ for $\tau_0\le\tau\le \tau_0+K, \; K={\rm const}>0$, and cross  the resonant surface
      ${\cal R}$ at some $\tau=\tau_*\in (\tau_0,  \tau_0+K), \; y=y_*, \; x=x_*$.

 \end{enumerate}

 Condition B) implies that the equation of the resonant surface can be presented in the form $I=a(y,x)$ with some smooth function $a(y,x)$. Then  $\omega_0(a(y,x), y,x)\equiv 0$.   Denote
\begin{eqnarray}
\alpha (y,x)&=& \frac{\p^2 H_0(a(y,x),y,x)}{\p I^2}, \label{alpha}\\
b(y,x)&=&\frac{\p a(y,x)}{\p x} \frac{\p H_0(a(y,x),y,x)}{\p y}-\frac{\p a(y,x)}{\p y} \frac{\p H_0(a(y,x),y,x)}{\p x}.\label{beta}
\end{eqnarray}
Condition B) implies that $\alpha(y,x)\ne0$.
Condition C) implies that $b(y,x)\ne0$.

For the sake of being definite assume  that  $\alpha(y,x)>0$, that $I>a(y,x)$ in the domain $\omega_0(I, y,x)>0$, and that $b(y,x)>0$. Then  condition C) implies that we should consider  the motion that starts at the domain $\omega_0>0$, i.e. $I_0>a(y_0,x_0)$.

\medskip

The following function is important in dynamics near the resonance:
\begin{equation}
\label{eq_F}
F(\vfi, y,x)=b(y,x)\vfi+H_1(a(y,x), \vfi, y,x, 0).
\end{equation}
We will assume that the following genericity condition is satisfied (its origin is explained in Section \ref{s_pendulum}).

\begin{enumerate}[label=\Alph*)]
  \setcounter{enumi}{4} 
  \item All critical points of the function $F$, considered as a function of $\vfi$, are non-degenerate: if $\frac{\p F}{\p \vfi}=0$ at some point, then  $\frac{\p^2 F}{\p \vfi^2}\ne 0$ at this point. Moreover, the values of $F$ at the critical points are different from each other.
\end{enumerate}
 Note that if $F$ does not have critical points, this condition is satisfied.

\begin{proposition}
\label{prop1}
Fix any number $r>1$. If the initial condition for a solution of system (\ref{e_system})  does not belong to some exceptional set ${\cal V}_r$ of  measure $O(\eps^r)$, then this solution crosses the resonan surface at some moment $\tau_e$ of the slow time and at a value  $\vfi_e$ of the phase. Moreover, $\tau_e =\tau_*+ O(\sqrt{\eps}\ln \nu)$, where  $\nu=\min (\frac 12, \nu_1+\eps)$, and $\nu_1>0$ is the minimum of distances of $\vfi_e$  from the points of local maximum  of the function $F(\vfi, y_*, x_*)$. If this function does not have critical points, then  ${\cal V}_r$  is empty and  $\tau_e =\tau_*+ O(\sqrt{\eps})$.
\end {proposition}
Here  $r>0$ can be taken arbitrary large, but the constant in the symbol $O(\sqrt{\eps}\ln \nu)$ depends on the choice of $r$.
Proposition \ref{prop1} is a direct corollary of Lemma \ref{lemma_after} formulated in Section \ref{upto}.

  \section {Formula for phase at the resonance}
  \label{s_main_formula}
  The second approximation of the averaging method is routinely used to approximately describe how the phase depends on the time far from resonances. It works on time intervals $\sim\frac 1\eps$ with accuracy $o(1)$ (see, e.g., \cite{bogol}). To apply the second-order averaging for the considered system, one should introduce new variables $(J,\psi, \eta, \xi)$  that are related to the old variables $(I,\vfi,y,x)$ via the canonical transformation of variables with the generating function
 \begin{equation}
 \label{gen_1}
 W (J,\vfi,\eta,x,\eps)= J\vfi+ \eps^{-1}\eta x+\eps  S_1(J,\vfi,\eta,x).
 \end {equation}
 Old and new variables are related as follows:
  \begin{eqnarray}
  \label{transform_1}
 I=J+\eps\frac {\p S_1(J,\vfi,\eta,x) }{\p\vfi}, \quad  \psi=\vfi+\eps\frac {\p S_1(J,\vfi,\eta,x) }{\p J},\\
 y=\eta+\eps^2\frac {\p S_1(J,\vfi,\eta,x) }{\p x}, \quad  \xi=x+\eps^2\frac {\p S_1(J,\vfi,\eta,x) }{\p\eta}. \nonumber
  \end{eqnarray}
Let
  $$
  {\cal H}_1(I,  y,  x)=\frac{1}{2\pi}\int_0^{2\pi}H_1(I,\vfi,y,x,0)d\vfi
  $$
  denote the average of the function $H_1$ over the fast phase at $\eps=0$.
The function $S_1$ is uniquely defined by the formula
  \begin{equation}
  \label{def_S_1}
  \frac{\p S_1(I,\vfi,y,x)}{\p\vfi}=\frac{1}{\omega_0(I,y,x)}\left[ {\cal H}_1(I,  y,  x)- H_1(I,\vfi,y,x,0)    \right]
  \end{equation}
  and the condition that the average of $S_1(I,\vfi,y,x)$ over $\vfi$ is equal to 0.

The dynamics of the variables $J,\psi, \eta, \xi$ is approximately described by the Hamiltonian system with Hamilton's function
  \begin{equation}
  {\cal H}_a(J, \eta,  \xi,\eps)= H_0( J,  \eta,  \xi)+\eps {\cal H}_{1}(J,  \eta, \xi).
  \end{equation}
   The corresponding differential equations are
   \begin{eqnarray}
   \label{i_system}
   \dot J&=&0,   \quad   \dot \eta=-\eps\frac{\p {\cal H}_a(J, \eta,  \xi,\eps) }{\p \xi},\quad   \dot \xi=\eps\frac{\p {\cal H}_a(J, \eta,  \xi,\eps) }{\p \eta},                                       \\
   \dot \psi&=&\frac{\p {\cal H}_a(J, \eta,  \xi,\eps) }{\p J}. \nonumber
   \end{eqnarray}
Thus, for the variables $\eta, \xi$ we have a Hamiltonian system with the Hamiltonian
$ \eps {\cal H}_a(J, \eta, \xi,\eps)$ that depends on the parameter $J={\rm const}$. Initial condition $(J_0, \eta_0, \xi_0)$  for this system should be obtained from the initial condition $(I_0,\vfi_0,  y_0, x_0)$ of the system (\ref{e_system}) using the transformation of variables (\ref{transform_1}). An approximation provided by these formulas is called {\it an improved adiabatic approximation}, the variable $J$ is {\it an improved adiabatic invariant}.

Denote
\begin{equation}
\omega_1(I,y,x)= \frac{\p  {\cal H}_1(I,  y,  x)}{\p I}.
\end{equation}
The frequency of the phase  $\psi$ in the system (\ref{i_system}) is $ \omega_0(J,  \eta, \xi)+\eps  \omega_1(J, \eta,\xi)$.

\medskip
 Denote by $\tilde H_1(I,\vfi,y,x)$ the purely periodic (i.e. with average 0) part of $H_1(I,\vfi, y, x,0)$. Then
 $\tilde H_1(I,\vfi,y,x)=H_1(I,\vfi, y, x,0)-{\cal H}_1(I,  y,  x)$. Denote
\begin{equation}
\label{eq_Xi}
\Xi=\frac{1}{2\pi}\left( \vfi_e+ \frac{\tilde H_1(I_0, \vfi_e, y_*,x_*)}{b(y_*,x_*)}     \right).
\end{equation}
Denote by $J_0, \eta_a(\tau), \xi_a(\tau)$ the solution of the differential equations \eqref{i_system} for  $J, \eta, \xi $  with initial condition $(J_0, \eta_0, \xi_0)$. Denote by $\tau_{*,a}$ the moment of the slow time when this solution crosses the resonant surface: $a(\eta_a(\tau_{*,a}), \xi_a(\tau_{*,a}))=J_0$.

\begin{theorem}
\label{theo_1}

If the initial point $(I_0,\vfi_0,  y_0, x_0)$  does not belong to the exceptional set ${\cal V}_r$, then the following asymptotic formula for $\Xi$ is valid:
\begin{equation}
\begin{aligned}
\label{Xi_a}
2\pi \Xi=&\vfi_0+\frac{1}{\eps}\int_{\tau_0}^{\tau_{*,a}} \left(    \omega_0(J_0,  \eta_a(\tau), \xi_a(\tau))+\eps  \omega_1(J_0,  \eta_a(\tau), \xi_a(\tau) \right)d\tau
+  O(\sqrt{\eps}\ln\nu).
\end{aligned}
\end{equation}
Here the objects  $r, {\cal V}_r, \nu$ are those introduced in Proposition \ref{prop1} and are associated with   points  of local maximum of the function  $F(\vfi, y_*, x_*)$. If $F$ does not have critical points, then ${\cal V}_r$  is  empty and  $\nu=1/2$, i.e. the error term in (\ref{Xi_a}) is
$ O(\sqrt{\eps})$.

\end{theorem}
This theorem is the main result of the paper. It is proven in Section \ref{main_proof} on the basis of some preliminary estimates and constructions contained in Sections \ref{s_pendulum}, \ref{preliminary}. Proofs of the preliminary estimates and auxiliary lemmas are given in Sections \ref{prelim_proof}, \ref{est_int}.

\medskip
{\bf Remark.}  Theorem \ref {theo_1} does not describe  initial conditions that belong to the exceptional set ${\cal V}_r$. However, this theorem determines  which initial conditions  are guaranteed to be outside ${\cal V}_r$. Let $\vfi_j, j=1,2, \ldots, k$, be points of maximum of the function $F(\vfi, y_*, x_*)$ on the interval $[0, 2\pi)$.  Denote $\tilde F$ the purely periodic part of the function $F$. Denote $\hat \Xi_j,\,   j=1,2, \ldots, k$,  fractional parts of values  $ \tilde F(\vfi_j, y_*, x_*)/(2\pi b(y_*, x_*))$.  Denote as  $\Xi_a$ the principal term  in expression for  $\Xi$ given by (\ref{Xi_a}). Let  $\hat \Xi_a$ be the fractional part of $\Xi_a$.   Then there exists a constant  $c_a>0$ such that if the initial point   $(I_0,\vfi_0,  y_0, x_0)$  is such that the value $\hat \Xi_a$ is at the distance bigger than  $c_a\sqrt{\eps}|\ln \eps| $  from each of values $\hat \Xi_j, j=1,2, \ldots, k$, then this initial point   does not belong to the exceptional set ${\cal V}_r$.

\medskip
Value $\Xi$ as a characteristic of a passage through resonance was introduced in \cite{Nei1997}.  It is determined modulo 1 (because the value of phase  $\vfi_e$ at the resonance is determined modulo $2 \pi$). We call $\Xi$ {\it pseudophase} in analogy with the terminology used in description of separatrices crossings \cite{NV_2005}.  The phase at the resonance $\vfi_e$ is uniquely determined by the pseudophase $\Xi$. If the function $F$ (\ref{eq_F}), considered as a function of $\vfi$,  does not have critical points, then one can solve  relation (\ref{eq_Xi}) for $\vfi_e$. Together with Theorem 1 this gives an asymptotic formula for $\vfi_e$ with an accuracy  $O(\sqrt{\eps}\,)$.
If  the function $F$ (\ref{eq_F}) has critical points, then $\vfi_e$ can't belong to some intervals  (this is explained in Section \ref{s_pendulum}). In this case relation (\ref{eq_Xi}) should be solved for $\vfi_e$ on some intervals of monotonicity of $F$ to get an asymptotic  formula for $\vfi_e$, see Section \ref{s_pendulum}.  The phase $\vfi_e$ will be  determined with an accuracy $O(\sqrt\eps\,)$ far from critical points of $F$.  The accuracy decays close to critical points of $F$. It is $O(\eps^{\frac 14}|\ln \eps|^{1/2})$ when $\Xi$ is $ \sqrt{\eps\,}|\ln \eps|$-close to a critical value $\tilde F_c$ of $\tilde F(\vfi, y_*, x_*)$ and is $O\left(\cfrac{\sqrt{\eps}\ln|\Xi-\tilde F_c|}{( \sqrt{|\Xi-\tilde F_c|}}\right)$ when  $|\Xi-\tilde F_c|>\sqrt{\eps}|\ln \eps|$.
(Logarithmic terms can be omitted for $\vfi_e$  near a critical point $\vfi_c$ of the function $F$  if the phase $\vfi$  is measured from this critical point, i. e., if   $\vfi_c\equiv 0$.)

\medskip
The problem of calculation of pseudophase $\Xi \ {\rm mod} \ 1$ (or phase $\vfi_e \ {\rm mod} \ 2 \pi$) for given initial conditions does not have a deterministic answer in the  limit as $\eps\to 0$. However, one can interpret $\hat \Xi={\rm Frac}\,(\Xi)$, the fractional part of $\Xi$,  as a random value and calculate a limiting probabilistic distribution of this value.  Such an approach to a problem with a small parameter was suggested   in \cite{LSN} and then was rediscovered and used by many authors (in particular, in \cite{Goldr}). A rigorous natural definition of what should be called a probability in such problems was suggested in \cite{Arn_63}. (There are also other natural definitions. In cases when results of calculation of probabilities using different definitions are known, these results coincide.)   According to \cite{Nei1997}, the value $\hat \Xi$
  has a uniform distribution on the interval $[0,1]$ in the limit as $\eps\to 0$. This result is obtained by calculation of phase fluxes, without knowledge of asymptotic formula for $\Xi$. Formula (\ref{Xi_a}) allows to give an alternative proof of this result (see Appendix).

\medskip

A general scheme of the proof of Theorem 1 is as follows. Inside $O(\sqrt{\eps})$-neighbourhood of the resonant surface one can use an expansion with respect to the distance from the resonant surface for description of the dynamics. It is known that the obtained ``pendulum-like" system has an approximate first integral (cf. Section \ref{s_pendulum}).  We introduce an auxiliary function  (this is the function $E$ in Section \ref{main_proof}) which is close to this first integral near the resonance. Due to this, one can get good estimates of change  of this function near the resonance: its change is
$O(\sqrt{\eps}\ln\nu)$ or $O(\sqrt{\eps})$ in $\sqrt{\eps}$-neighbourhood of the resonant surface. Outside $\sim\sqrt{\eps}$-neighbourhood of the resonant surface one can make a transformation of variables which excludes fast phase from the Hamiltonian  up to a prescribed order in $\eps$. Using such transformations is a standard approach in the averaging method \cite{bogol}. This allows to calculate the change of function $E$  outside $\sqrt{\eps}$-neighbourhood of the resonance with an accuracy $O(\sqrt{\eps})$. Combining the estimates near and far from the resonance we get an asymptotic formula for change of $E$ from the initial value to arrival into resonance with an accuracy $O(\sqrt{\eps}\ln\nu)$ or $O(\sqrt{\eps})$. This provides an equality that includes value
$\Xi$ (\ref{eq_Xi}). Solving this equality for $\Xi$ gives relation (\ref{Xi_a}) in Theorem 1.

\section{Expansion near the resonant surface}
\label{s_pendulum}
 Study of dynamics close to resonance can be reduced to study of behaviour of a perturbed pendulum-like system. Such a reduction was used in a number of papers starting from \cite{Chir1959}; see, e.g., \cite{AKN} and references therein. We will perform this reduction in a Hamiltonian form following \cite{Nei_steklov}.

Expansion of Hamiltonian (\ref{H_exact}) over $(I-a(y, x))$ is
 \begin{eqnarray}
   H
   =H_0(a( y, x),  y, x)+\frac12 \alpha( y, x) (I-a( y, x))^2
+\eps  H_{1}( a( y, x),\vfi,  y, x, 0) \\
+ O((I-a( y, x))^3)+\eps O(I-a( y, x))+O(\eps^2),\phantom{***********}\nonumber
\end{eqnarray}
where the function $\alpha$  is given by formula (\ref{alpha}).
Make a canonical transformation of variables $(I,\vfi, y, x)\mapsto(\sqrt{\eps}P, \gamma,\tilde y, \tilde x)$ with the generating function
\begin{equation}
 \label{gen_W}
   W_{res}(P, \vfi, \tilde y, x, \eps)=\tilde y\eps^{-1}x+(\sqrt{\eps}P+a(\tilde y, x))\vfi .
\end{equation}
Then $\gamma=\vfi, \  I=\sqrt{\eps}P+a(\tilde y,\tilde x)+O(\eps)$. Other old variables are expressed via new ones as follows:
\begin{eqnarray}
y=\tilde y+\eps\dfrac{\p a(\tilde y,\tilde x)}{\p\tilde x}\vfi+O(\eps^2), \quad
      x=\tilde x-\eps\dfrac{\p a(\tilde y,\tilde x)}{\p\tilde y}\vfi+O(\eps^2).
\end{eqnarray}
Note that in these estimates we  assume that $\vfi=O(1)$.
The Hamiltonian (\ref{H_exact}) in the new variables takes the form
\begin{eqnarray}
\label{H_expanded}
   H&=&H_0(a(\tilde y, \tilde x),\tilde y,  \tilde x)\\
   &+&\eps\left(  \frac 12\alpha(\tilde y, \tilde x)P^2+b(\tilde y, \tilde x)\vfi+H_1(a(\tilde y, \tilde x),\vfi,\tilde y, \tilde x,0)    \right)\nonumber\\
   &+&\eps^{\frac32}O(P^3)+\eps^{\frac32} O(P)+O(\eps^2),\nonumber
\end{eqnarray}
where the function $b$  is given by formula (\ref{beta}). Consider  motion in a strip near the resonant surface where $P=O(1)$, and neglect terms $O(\eps^{\frac32})$. Introduce the new time $\theta=\sqrt{\eps} t$. For dynamics of $\tilde y, \tilde x$ considered as canonically conjugate variables we get  a Hamiltonian system with  Hamilton's function $ \sqrt{\eps}\Lambda(\tilde y, \tilde x)$, where $\Lambda=H_0(a(\tilde y, \tilde x),\tilde y,  \tilde x)$.
Dynamics of $P, \vfi$ in the principal approximation can be considered for frozen $\tilde y, \tilde x$. Then for $P, \vfi$  we get the Hamiltonian system with  Hamilton's function
\begin{equation}
\label{def_E}
{\cal E}(P, \vfi, \tilde y, \tilde x)=\frac 12\alpha(\tilde y, \tilde x)P^2+b(\tilde y, \tilde x)\vfi+H_1(a(\tilde y, \tilde x),\vfi,\tilde y, \tilde x,0),
\end{equation}
i.e.
\begin{equation}
\label{def_E1}
{\cal E}(P, \vfi, \tilde y, \tilde x)=\frac 12\alpha(\tilde y, \tilde x)P^2+F(\vfi, \tilde y, \tilde x),
\end{equation}
where the function $F$  is given by formula (\ref{eq_F}). Here $P, \vfi$ are canonically conjugate variables (note that in (\ref {H_expanded}) canonically conjugate variables are  $\sqrt{\eps}P, \vfi$). 

\begin{figure}
\begin{center}
            \includegraphics[scale=0.35]{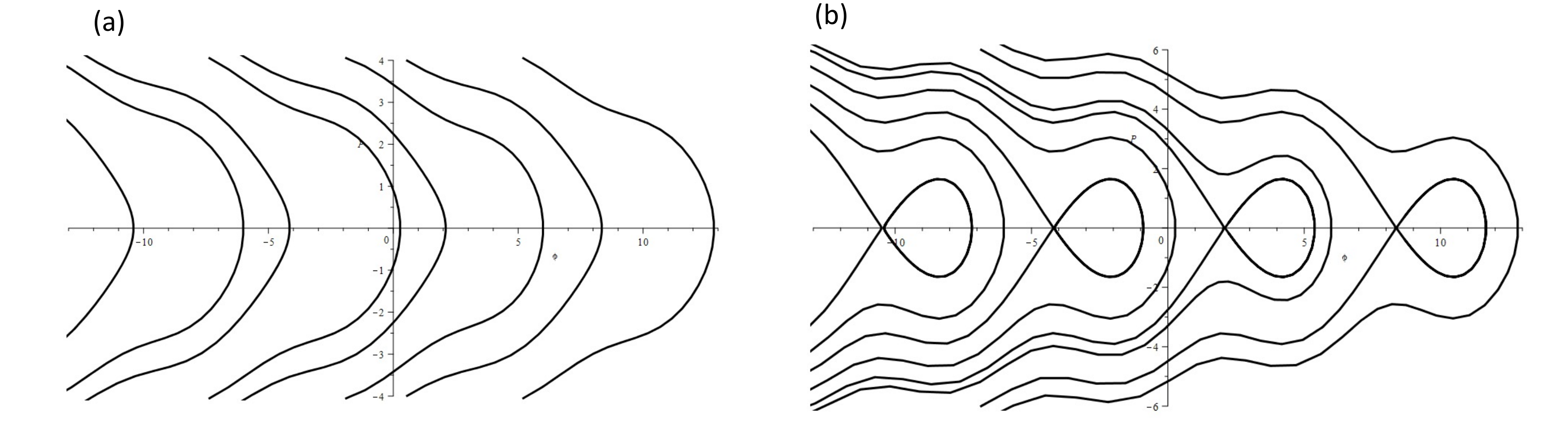}
  \end{center}
  \vskip -0.4cm
           \caption{Phase portraits of pendulum-like systems: (a) without oscillatory domains, (b)  with oscillatory domains.}
 \label{portraits}
\end{figure}

Hamiltonian (\ref{def_E}) describes a one-dimensional motion in  $2\pi$-periodic potential $H_1(a(\tilde y, \tilde x),\vfi,\tilde y, \tilde x,0)$ under the action of the constant torque $b(\tilde y, \tilde x)$ (a ``pendulum-like system'').     There are two types of  phase portraits of this Hamiltonian:  without oscillatory domains (Fig. \ref{portraits}, a) and with oscillatory domains (Fig.  \ref{portraits}, b).
Condition E) in Section \ref{s_assumptions} implies that unstable equilibria  in Fig.  \ref{portraits}, b are non-degenerate saddles, and that separatrices do not connect different saddles \footnote{Condition E) can be weakened: it is enough to consider only critical points that correspond to saddles at boundaries of oscillatory domains.}.

\medskip
If  the phase portrait of the Hamiltonian ${\cal E}$ does not have oscillatory domains for $x=x_*, y=y_*$, then the phase at the resonance $\vfi_e$ can have any value   modulo $ 2\pi$.  To determine an approximate value of $\vfi_e$ it is enough to solve the relation (\ref{eq_Xi}) for $\vfi_e$ and use the approximate value for $\Xi$ given by formula (\ref{Xi_a}). If there are oscillatory domains, then, in the approximation (\ref{def_E1}),   values of  $\vfi_e$  can't correspond to points  in  oscillatory domains.  To find  $\vfi_e$ ${\rm mod}\ 2\pi$, one can solve relation (\ref{eq_Xi}) on intervals  of monotonicity of the function $F$ for $\vfi$ outside oscillatory domains and use the approximate value for $\Xi$ from (\ref{Xi_a}). However, due to the perturbation in (\ref{H_expanded}), the point $P=0, \varphi=\varphi_e$ can be in an oscillatory domain
 with value of $\Xi$ which is  $O(\sqrt\eps)$-close to the value of $F(\vfi, y_*,x_*)$  at a saddle point  on the boundary of this domain.

\section{Preliminary estimates}
\label{preliminary}
Fix an initial condition $(I_0, \vfi_0, y_0, x_0)$ at $t=t_0$ and the solution
$(I(t),\vfi(t), y(t),x(t))$  of  system (\ref {e_system}) with this initial condition. Rewrite these initial conditions and  solution in variables $(J,\psi, \eta, \xi)$ introduced by the transformation
(\ref {transform_1}). We get an initial condition $(J_0, \psi_0, \eta_0, \xi_0)$ and solution  $(J(t), \psi(t), \eta(t), \xi(t))$.  Consider solution $(J_0, \psi_a(\tau), \eta_a(\tau),\xi_a(\tau)) $     of  system (\ref {i_system}) with initial condition $(J_0, \psi_0, \eta_0, \xi_0)$  at $\tau=\tau_0 =\eps t_0$. Denote $\omega_a(\tau)=\omega_0(J_0, \eta_a(\tau), \xi_a(\tau))$.

\subsection{Estimates up to distance $\sim \sqrt{\eps}$ from the resonance}

\begin{lemma}
\label{lemma_before}
There exists a constant $C_1>0$ such that while $\omega_a(\tau)\ge C_1\sqrt{\eps}$ the solution $(J(t), \psi(t), \eta(t), \xi(t))$ is well defined  and the following estimates are satisfied:
\begin{eqnarray}
|J(t)-J_0|=O\left(\frac{\eps^2}{\omega^3_a(\tau)}\right),\quad  |\psi(t)-\psi_a(\tau)|=O\left(\frac{\eps}{\omega^2_a(\tau)}\right),\label{line_1_6.1}\\
|\eta(t)-\eta_a(\tau)|=O\left(\frac{\eps^2}{\omega^2_a(\tau)}\right),\quad  |\xi(t)-\xi_a(\tau)|=O\left(\frac{\eps^2}{\omega^2_a(\tau)}\right),
\label{line_2_6.1}\\
| I(t)-J(t)|=O\left(\frac{\eps}{\omega_a(\tau)}\right),\quad  |\vfi(t)-\psi(t)|=O\left(\frac{\eps}{\omega^2_a(\tau)}\right), \label{line_3_6.1}\\
| y(t)-\eta(t)|=O\left(\frac{\eps^2}{\omega^2_a(\tau)}\right),\quad  |x(t)-\xi(t)|=O\left(\frac{\eps^2}{\omega^2_a(\tau)}\right), \label{line_4_6.1}\\
\omega_0(J(t), \eta(t), \xi(t))>\frac{1}{2}\omega_a(\tau),\quad \omega_0(I(t), y(t), x(t))>\frac{1}{2}\omega_a(\tau). \label{line_5_6.1}
\end{eqnarray}

\end{lemma}

This lemma is proven in Section \ref{proof_lemma_before}.
\subsection{Estimates for slow variables up to arrival to  resonance}
\label{upto}
\begin{lemma}
\label{lemma_after}
(a) There exist a constant $C_2>C_1$  and a moment of time  $t_{N}$   such that
$$\omega_0(I(t_{N}), y(t_{N}), x(t_{N}))=C_2\sqrt {\eps},\ \omega_a(\eps t_{N})>2C_1\sqrt {\eps},\ \eps t_{N}=\tau_*+O(\sqrt{\eps}).
$$
(b)  Fix a number $r>1$. If the initial condition  $(I_0, \vfi_0, y_0, x_0)$ does not belong to an exceptional set ${\cal V}_r$ of measure $O(\eps^r)$, then the solution $(I(t), \vfi (t), y(t), x(t))$ arrives to the resonant surface at some $t=t_e$ (i.e.
$\omega_0(I(t_e), y(t_e), x(t_e))=0$, or, which is the same, $a(y(t_e), x(t_e))=I(t_e)$). Denote $\nu=\min (\frac 12, \nu_1+\eps)$, where $\nu_1>0$ is the minimum of distances of $\vfi(t_e)$ from the  points of maximum of  the function $F(\vfi, y_*, x_*)$.  Then
$\eps t_e=\tau_* +O(\sqrt{\eps}\ln\nu)$.  For $t_{N}  \le t\le t_e$
the following estimates are satisfied:
\begin{equation}
\begin{aligned}
\label{est_br_1}
I(t)&=I_0+O(\sqrt{\eps}\ln\nu),
\\ y(t)&=\eta_a(\tau)+O(\eps\ln^2\nu), \; x(t)=\xi_a(\tau)+O(\eps\ln^2\nu).
 \end{aligned}
\end{equation}
(c) If the function $F$  does not have critical points, then for any initial conditions the solution arrives to the resonance surface, and one can omit $\ln \nu$ in the previous estimates, i.e.
\begin{equation}
\begin{aligned}
\label{est_br_2}
\eps t_e&=\tau_* +O(\sqrt{\eps}), \ I(t)=I_0+O(\sqrt{\eps}), \\
y(t)&=\eta_a(\eps t)+O(\eps), \  x(t)=\xi_a(\eps t)+O(\eps).
\end{aligned}
\end{equation}
\end{lemma}

This lemma is proven in Section \ref{proof_lemma_after}.

\section{Proof of Theorem 1}
\label{main_proof}
We consider solution $(I(t),\vfi(t), y(t),x(t))$  of the system (\ref {e_system}). We assume that initial condition  $(I_0, \vfi_0, y_0, x_0)$ at $t=t_0$ of this solution does not belong to the exceptional set ${\cal V}_r$. According to Proposition \ref{prop1}, this solution crosses the resonant surface at some moment of time $t_e$. We denote $I_e= I(t_e), \vfi_e= \vfi(t_e), y_e=y(t_e),x_e=x(t_e),
\tau_0 =\eps t_0, \tau_e =\eps t_e$. Crossing the resonance means that $\omega_0(I_e,y_e,x_e)=0$ or, equivalently, $a(y_e, x_e)=I_e$.   Without loss of generality we put $t_0=0, \tau_0=0$.  We use prime to denote derivative with respect to $\tau$.

We use notation of previous sections: $(J(t), \psi(t), \eta(t), \xi(t))$ is the solution represented  in variables (\ref {transform_1}), $(J_0, \psi_0, \eta_0, \xi_0)$ is its initial condition, $(J_0, \psi_a(\tau), \eta_a(\tau),\xi_a(\tau)) $ is the solution of system (\ref {i_system}) with the initial condition $(J_0, \psi_0, \eta_0, \xi_0)$  at  $\tau=0$,  and $\omega_a(\tau)=\omega_0(J_0, \eta_a(\tau), \xi_a(\tau))$. There is a moment of the slow time $\tau=\tau_{*,a}$ such that $\omega_a(\tau_{*,a})=0$.  Proposition \ref{prop1} implies that $\tau_e= \tau_{*,a}+O(\sqrt{\eps}\ln \nu)$.

The expansion of $H_0(I,y,x)$ over $(I-a(y, x))$ is
\begin{equation}
\label{exp_H_0}
H_0(I,y,x)=\Lambda(y,x)+\frac 12\alpha(y,x)(I-a(y,x))^2+M(I,y,x),
\end{equation}
where $M(I,y,x)=O((I-a(y,x))^3)$.
 From (\ref{e_system}) we get
\begin{equation}
\dot \vfi= \alpha(y,x)(I-a(y,x))+\frac{\p M(I, y,x)}{\p I}+\eps\frac{\p H_1(I,\vfi, y,x,0)}{\p I}+O(\eps^2).
\end{equation}

Denote
\begin{equation*}
\begin{aligned}
B_0(I, y,x)=&\frac{\p a(y,x)}{\p x}\frac{\p H_0(I,y,x)}{\p y}-\frac{\p a(y,x)}{\p y}\frac{\p H_0(I,y,x)}{\p x}\\
=&b(y,x)+\frac{\p a(y,x)}{\p x}\frac{\p N(I,y,x)}{\p y}-\frac{\p a(y,x)}{\p y}\frac{\p N(I,y,x)}{\p x},\\
B_1(I,\vfi ,y,x)=&\frac{\p a(y,x)}{\p x}\frac{\p H_1(I,\vfi,y,x,0)}{\p y}-\frac{\p a(y,x)}{\p y}\frac{\p H_1(I,\vfi, y,x,0)}{\p x},\\
 b_{*,a}&=b(\eta_a(\tau_{*,a}), \xi_a(\tau_{*,a})), \ N(I,y,x)=O((I-a(y,x))^2).
 \end{aligned}
\end{equation*}

\medskip

Consider the function:
\begin{equation}
E=\frac{1}{2\eps}\alpha(y,x)(I-a(y,x))^2+H_1(I,\vfi, y,x,\eps)+b_{*,a}\vfi.
\end{equation}
This function is an approximate first integral of motion in system (\ref {e_system})  near the resonant surface. Calculate derivative $\dot E$ of this function along solution  $(I(t),\vfi(t), y(t),x(t))$  and take the integral of this derivative from 0  to $t_e$.  Compare two expressions of this integral.
On the one hand,
\begin{equation}
\begin{aligned}
 \label{int_dot_E_1}
&\int^{t_e}_0\dot E dt=E|_{t=t_e}-E|_{t=0}\\
=&H_1(I_e,\vfi_e,y_e,x_e,\eps)+b_{*,a}\vfi_e-\frac1{2\eps}\alpha(y_0,x_0)(I_0-a(y_0,x_0))^2\\
&-H_1(I_0,\vfi_0,y_0,x_0,\eps)-b_{*,a}\vfi_0.
\end{aligned}
\end{equation}
For what follows we need to express $I_0$ via $J_0$
:
$$
I_0=J_0-\eps\frac{\tilde H_1(I_0,\vfi_0,y_0,x_0)}{\omega_0(I_0,y_0,x_0)}+O(\eps^2).
$$

Thus we have
\begin{equation}
\begin{aligned}
 \label{int_dot_E_1_1}
&\int^{t_e}_0\dot E dt=H_1(I_e,\vfi_e,y_e,x_e,0)+b_{*,a}\vfi_e-\frac1{2\eps}\alpha(y_0,x_0)(J_0-a(y_0,x_0))^2\\
&+\alpha(y_0,x_0) (J_0-a(y_0,x_0))\frac{\tilde H_1(I_0,\vfi_0,y_0,x_0)}{\omega_0(I_0,y_0,x_0)}\\
 & -H_1(I_0,\vfi_0,y_0,x_0,\eps)-b_{*,a}\vfi_0 +O(\eps)\\
 =&H_1(I_e,\vfi_e,y_e,x_e,0)+b_{*,a}\vfi_e-\frac1{2\eps}\alpha(y_0,x_0)(J_0-a(y_0,x_0))^2\\
&-\frac{\p M}{\p I}(I_0,y_0,x_0)\frac{\tilde H_1(I_0,\vfi_0,y_0,x_0)}{\omega_0(I_0,y_0,x_0)}
 -{\cal H}_1(I_0,y_0,x_0) -b_{*,a}\vfi_0 +O(\eps).
\end{aligned}
\end{equation}

On the other hand, we should integrate  from 0 to $t_e$ the expression
\begin{equation}
\begin{aligned}
\dot E=&\frac 1\eps\alpha(y,x)(I-a(y,x))(\dot I-\eps a'(y,x))+\frac 12\alpha '(y,x)(I-a(y,x))^2\\
&+\frac{\p H_1}{\p I}\dot I+\frac{\p H_1}{\p \vfi}\dot\vfi+\frac{\p H_1}{\p y}\dot y+\frac{\p H_1}{\p x}\dot x+b_{*,a}\dot\vfi
.
\end{aligned}
\end{equation}
Here ``prime'' denotes the derivative with respect to the slow time $\tau$.

Substitute expressions for $\dot I, \dot \vfi, \dot y, \dot x$. Some terms are cancelled, which is expected because $E$ is the first integral of the expanded near the resonance system for frozen slow variables. We finally get
\begin{equation}
\begin{aligned}
 \label{dot_E}
\dot E=&-\alpha(y,x)(I-a(y,x)) (b(y,x)-b_{*,a})-\eps \alpha(y,x)(I-a(y,x))B_1(I,\vfi,y,x)\\
-&\alpha(y,x)(I-a(y,x))\left(\frac{\p a(y,x)}{\p x}\frac{\p N(I,y,x)}{\p y}
-\frac{\p a(y,x)}{\p y}\frac{\p N(I,y,x)}{\p x}\right)\\
+&\frac 12(I-a(y,x))^2\left(\frac{\p \alpha(y,x)}{\p x}\frac{\p (H_0+\eps H_1)}{\p y}-\frac{\p \alpha(y,x)}{\p y}\frac{(\p H_0+\eps H_1)}{\p x}\right)\\
+&\frac{\p H_1}{\p \vfi}\frac{\p M}{\p I}
+\eps\left(-\frac{\p H_1}{\p y}\frac{\p H_0}{\p x}+\frac{\p H_1}{\p x}\frac{\p H_0}{\p y}\right)
+b_{*,a}\left(\frac{\p M}{\p I}+\eps\frac{\p H_1}{\p I}\right)+O(\eps^2).
\end{aligned}
\end{equation}
We should calculate
\begin{equation*}
\int^{t_e}_0\dot E dt=\int_0^{t_N}\dot E dt+\int^{t_e}_{t_N}\dot E dt,
\end{equation*}
where $t_N$ is the moment of time introduced in Lemma \ref{lemma_after}.
\begin{lemma}
\label{last_integral}
\begin{equation*}
\int^{t_e}_{t_N}\dot E dt= O(\sqrt{\eps}\ln \nu).
\end{equation*}
\end{lemma}
The proof is in Section \ref{est_int}.

Our goal now is to estimate the integral of $\dot E$ from 0 to $t_N$.  On this time interval
we have $\omega_a(\tau)\ge C_1\sqrt{\eps}$. Hence, estimates of Lemma \ref{lemma_before} are valid. Using   (\ref{transform_1}) we get
\begin{equation}
\begin{aligned}
\label{Part8_con_1}
I =&J_0+\eps\frac {\p S_1}{\p\vfi}\left(J_0,\psi,\eta_a,\xi_a\right) +O\left(\frac{\eps^2}{\omega^3_a(\tau)}\right), \\
\vfi= &\psi-\eps\frac {\p S_1}{\p I}\left(J_0,\psi,\eta_a, \xi_a\right) +O\left(\frac{\eps^2}{\omega^4_a(\tau)}\right),\\
y=&\eta_a+O\left(\frac{\eps^2}{\omega^2_a(\tau)}\right),\ x=\xi_a+O\left(\frac{\eps^2}{\omega^2_a(\tau)}\right).
\end{aligned}
\end{equation}
We will estimate integrals of terms in  (\ref{dot_E}) with the help of lemmas formulated below. They are proven in Section \ref{est_int}. In these lemmas, if some function of $I,\vfi, y,x,\psi$ is under an  integral over time $t$, then $I=I(t),\vfi=\vfi(t), y=y(t),x=x(t), \psi=\psi(t)$. \begin{lemma}
\label{int_with S1}
Let a function $L$ have a form
$$L(J_0,\psi,\eta_a,\xi_a)=(J_0-a(\eta_a,\xi_a))
\left(\frac {\p S_1}{\p\vfi}\left(J_0,\psi,\eta_a,\xi_a\right)
\right)L_1\left(J_0,\psi,\eta_a,\xi_a\right)
$$
 where $L_1$ is a smooth function. Then
\begin{equation*}
\eps\int_{0}^{t_N}L(J_0,\psi,\eta_a,\xi_a) dt=O(\sqrt \eps ).
\end{equation*}
\end{lemma}

\begin{lemma}
\label{lemma_T0}
Let a function $L$ have a form  $L(I,y,x)=(I-a(y,x))^2L_1(I,y,x)$, where $L_1$ is a smooth function. Then
\begin{equation*}
\int_{0}^{t_N}L(I,y,x)dt=\frac{1}{\eps}\int_0^{\tau_{*,a}}L(J_0, \eta_a, \xi_a)d\tau +O(\sqrt \eps ).
\end{equation*}
\end{lemma}

\begin{lemma}
\label{lemma_T1}
Let ${L}(I,y,x,\eps)$ be a smooth  function.
Then
$$
\eps\int_{0}^{t_N}L(I,y,x,\eps)dt=\int_{0}^{\tau_{*,a}}{L}(J_0,\eta_a,\xi_a,0)d\tau +O(\sqrt{\eps}).
$$
\end{lemma}

\begin{lemma}
\label{lemma_T2}
Let ${L}(I,\vfi,y,x, \eps)$ be a smooth $2\pi$-periodic function of $\vfi$ with average 0.
Then
$$
\eps\int_{0}^{t_N}L(I,\vfi,y,x, \eps)dt=O(\sqrt{\eps}).
$$
\end{lemma}

\begin{lemma}
\label{lemma_T3}
Let
$$
L(I,y,x)=\alpha (y,x)(I-a(y,x)) (b(y,x)-b_{*,a}).
$$
Then
\begin{equation*}
\begin{aligned}
\int_{0}^{t_N}L(I,y,x)dt= \frac{1}{\eps}\int_{0}^{\tau_{*,a}}L(J_0,\eta_a,\xi_a)dt
+O(\sqrt \eps).
\end{aligned}
\end{equation*}
\end{lemma}

\begin{lemma}
\label{lemma_T4}
Let  $K=K(I,y,x)$ be a smooth function such that $K(I,y,x)= O((I-a(y,x))^2)$. Then
\begin{equation*}
\begin{aligned}
\int^{t_N}_0 \frac{\p H_1}{\p \vfi}(I,\vfi,y,x,\eps) K(I,y,x)dt
=-\frac{\tilde H_1(I_0, \vfi_0,y_0,x_0)}{\omega_0(I_0, y_0, x_0)} K(I_0,y_0,x_0)+O({\sqrt \eps}).
\end{aligned}
\end{equation*}
\end{lemma}
With these lemmas we can estimate integrals of all terms in (\ref{dot_E}). According to Lemma \ref{lemma_T4}, we have
\begin{equation}
\begin{aligned}
\label{int_dot_E_sep_from_line 4}
\int^{t_N}_0 \frac{\p H_1}{\p \vfi}\frac{\p M}{\p I} dt=
 -\frac{\tilde H_1(I_0, \vfi_0,y_0,x_0)}{\omega_0(I_0, y_0, x_0)} \frac{\p M}{\p I}(I_0,y_0,x_0)+O(\sqrt\eps).
\end{aligned}
\end{equation}
For all other terms in (\ref{dot_E}) one should replace $I,y,x$ with $J_0, \eta_a, \xi_a$  and  average over $\vfi$. The integral should be taken from 0 to $\tau_{*,a}$. The accuracy of this approximation is $O(\sqrt\eps\ln \nu)$. This means that for approximate calculation of change of $E$ one can replace $I,y,x,\vfi$ and $H_1$ with  $J_0, \eta_a, \xi_a, \psi_a$ and ${\cal H}_1$, but an additional term (\ref {int_dot_E_sep_from_line 4}) appears. (This can be checked also by a direct substitution of terms.) Thus we have
\begin{equation}
\begin{aligned}
\label{second_integral_form}
&E|_{t=t_e}-E|_{t=0}=-\frac{\tilde H_1(I_0, \vfi_0,y_0,x_0)}{\omega_0(I_0, y_0, x_0)} \frac{\p M}{\p I}(I_0,y_0,x_0)\\
& -\frac{1}{2}\alpha(y_0,x_0)(J_0-a(y_0,x_0))^2
+{\cal H}_1(I_0,y_*,x_*)-{\cal H}_1(I_0,y_0,x_0)\\
+&\frac{ b_{*,a}}{\eps}\int^{\tau_{*,a}}_0\left(  \omega_0(J_0,\eta_a,\xi_a)  +\eps \omega_1(J_0,\eta_a,\xi_a)    \right)d\tau
+O(\sqrt\eps \ln \nu).
\end{aligned}
\end{equation}
Now we can equate two expression for $E|_{t=t_e}-E|_{t=0}$, those in (\ref{int_dot_E_1_1}) and in
(\ref{second_integral_form}). Several terms are cancelled, and we get
\begin{equation*}
\begin{aligned}
&b_{*,a}\vfi_e -b_{*,a}\vfi_0+H_1(I_e,\vfi_e,y_e,x_e,0)
={\cal H}_1(I_0,y_*,x_*)\\
+&\frac{ b_{*,a}}{\eps}\int^{\tau_{*,a}}_0\left(  \omega_0(J_0,\eta_a,\xi_a)  +\eps \omega_1(J_0,\eta_a,\xi_a)    \right)d\tau
 +O(\sqrt\eps \ln \nu).
\end{aligned}
\end{equation*}
This implies
\begin{equation}
\begin{aligned}
\label{first_part}
&\vfi_e +\frac{\tilde H_1(I_0,\vfi_e,y_*,x_*)
 }{b_{*}}
 =\vfi_0+\frac{ 1}{\eps}\int^{\tau_{*,a}}_0\left(  \omega_0(J_0,\eta_a,\xi_a)
+\eps \omega_1(J_0,\eta_a,\xi_a)    \right)d\tau\\
&+O(\sqrt\eps\ln \nu).
\end{aligned}
\end{equation}
According to the definition of pseudophase $\Xi$ (\ref{eq_Xi}), this gives the assertion of Theorem  \ref{theo_1}.

\section{Numerical test }
We performed a numerical check of accuracy of the formula for pseudophase  at the resonance. This was done for the Hamiltonian
\begin{equation}
\begin{aligned}
H=y + \frac 12 (I-x^2)^2+\eps(\frac 12+I)\sin \vfi
\end{aligned}
\end{equation}
(pairs of canonically conjugate variables are $I, \vfi$ and $y, \eps^{-1}x$).
Differential equations of motion for variables $I, \vfi, x$ are
\begin{equation}
\begin{aligned}
\label{exper_example}
 \dot I= -\eps(\frac 12+I) \cos\vfi, \  \dot \vfi = (I-x^2)+\eps \sin\vfi,  \  \dot x  = \eps.
\end{aligned}
\end{equation}
Equations for $I, \vfi$ form a Hamiltonian system that depends on the slow time $x \equiv \tau=\eps t$.   In notation of Sections \ref{s_outline}, \ref{s_main_formula} we have
\begin{equation}
\omega_0=I-\tau^2, \quad \omega_1=0,\quad  H_1=(\frac 12+I)\sin \vfi,\quad  {\cal H}_1=0, \ b=2\tau \, .
\end{equation}
Variables $I$ and $J$ are related as (cf. (\ref{transform_1}), (\ref{def_S_1}))
$$
I=J -\eps\frac{(\frac 12+J)\sin \vfi}{J-\tau^2} \,.
$$
The resonant surface has the equation $I=\tau^2$. We consider motion that starts at $t=0$ with $I=I_0, \vfi=\vfi_0$. The moment of arrival to the resonance in the improved adiabatic approximation is  $\tau_{*,a}=\sqrt{J_0}$ , where $J_0$ is the value of $J$ at $\tau =0$. The Hamiltonian of pendulum-like system is (cf. (\ref{def_E}))
$$
{\cal E}=\frac{1}{2}P^2+(\frac 12+\tau^2)\sin \vfi +2\tau\, \vfi.
$$
In what follows we have $\tau=\tau_{*,a}=1$ in this Hamiltonian. There are no equilibria in its phase portrait Fig.  \ref{portraits}, a. The error  term in Theorem \ref{theo_1} is $O(\sqrt{\eps})$.

\medskip
 The definition of the pseudophase  $\Xi$  (\ref{eq_Xi}) gives
\begin{equation}
\label{Xi_def_e}
\Xi=\frac{1}{2\pi}\left( \vfi_e + \frac{1}{2\tau_*}(\frac12+\tau_*^2)\sin \vfi  \right)\,.
\end {equation}
The main term $\Xi_{theor}$ of the asymptotic formula (\ref{Xi_a})  for calculation of value of pseudophase  at arrival to resonance reads as
$$
\Xi_{theor}=\frac{1}{2\pi}\left( \vfi_0+\frac{1}{\eps}\int_0^{\tau_{*,a}}(J_0-\tau^2)d \tau\right).
$$
Thus, it can be calculated explicitly:
\begin{equation}
\label{Xi_t}
 \Xi_{theor}=\frac{1}{2\pi}\left(\vfi_0+\frac{2J_0^{3/2}}{3\eps} \right).
\end{equation}
For the numerical tests we used  $J_0=1$,  values of the small parameter  $\eps= 0.001\times0.5^k, k=1, \ldots,10$,  and  250 values  of the initial phase  $\vfi_0$ evenly distributed between $0$ and $2\pi$. Values  $\Xi$ (\ref{Xi_def_e}) were obtained by means of solving equations (\ref {exper_example}) numerically.  The numerical values were compared with  $\Xi_{theor}$ given by  formula (\ref{Xi_t}).
\begin{table}[!ht]
\centering
\begin{tabular}{c c c c c c }
\toprule
$\eps$ & $0.0005$& $0.001\times0.5^2$& $0.001\times0.5^3$& $0.001\times0.5^4$&$0.001\times0.5^5$ \\ [0.5ex]
\midrule
RMSE&0.0013758&0.0009735&0.0007024 &0.0004924 &0.0003457 \\[1ex]
\bottomrule
\end{tabular}
\begin{tabular}{c c c c c c}
\toprule
$\eps$ & $0.001\times0.5^6$&$0.001\times0.5^7$&$0.001\times0.5^8$&$0.001\times0.5^9$&$0.001\times0.5^{10}$\\ [0.5ex]
\midrule
RMSE&0.0002458 &0.0001735 &0.0001227&0.0000876&0.0000621  \\[1ex]
\bottomrule
\end{tabular}
\caption{Values of $\eps$ and RMSE}
\label{table:1}
\end{table}
For each $\eps$ we calculated the root mean squared error (RMSE) between these values over 250 values of the initial phase. To check the scaling of this error with  $\eps$ we   took logarithms of RMSE and $\eps$ and fit them with a  linear equation. The result is that $\log({\rm RMSE})=0.4977231\log(\eps)-2.8018347 $, and the fit image is shown below. This scaling is in a good agreement with the estimate $O(\sqrt{\eps})$ of the error term in  Theorem \ref{theo_1}.

 \begin{figure}
\begin{center}
            \includegraphics[scale=0.4, angle=0.0]{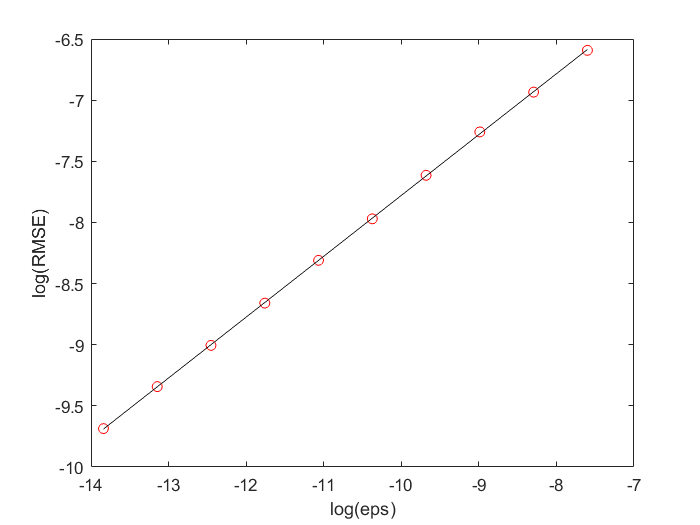}
  \end{center}
  \vskip -0.8cm
           \caption{Fit image}
 \label{portraits_fit}
\end{figure}

\medskip
\section{Proofs of preliminary estimates}
\label{prelim_proof}
\subsection{Proof of Lemma \ref{lemma_before}}
\label{proof_lemma_before}
The canonical perturbation theory (see, e.g., \cite{AKN}) provides a construction of  canonical transformations of variables such that the Hamiltonian in the new variables does not depend on the fast phase up to any prescribed order in $\eps$. We will use two such   transformations. They eliminate dependence on a fast phase up  to the first and the second order in $\eps$, respectively.

The dependence on a fast phase up  to the first order in $\eps$ is eliminated by the transformation to the variables  $(J,\psi, \eta, \xi)$ introduced in Section \ref{s_main_formula}. This transformation is determined by the generating function (\ref {gen_1}), function $S_1$ is defined by formula (\ref{def_S_1}) and by the condition that the average of $S_1$ over $\vfi$ is 0.  New and old variables are related via formulas (\ref{transform_1}). The Hamiltonian in the new variables has the form
\begin{equation}
      \label{H_transformed_S1}
     {{\cal H}}(J, \psi, \eta, \xi, \eps)
     = H_0( J,\eta, \xi)+ \eps {\cal H}_1(J, \eta, \xi)
      + \eps^2 {\cal H}_{2}( J, \psi,\eta, \xi,\eps).
\end{equation}

In order to describe behaviour of the variable $J$ with a better accuracy we need  also the transformation of variables that eliminates dependence on a fast phase up  to the second order in $\eps$.
The new variables $(\hat J,\hat\psi,\hat \eta,\hat \xi)$ are related with the old variables $(I,\vfi,y,x)$ via the canonical transformation of variables with a generating function
 \begin{equation}
 \hat W (\hat J,\vfi,\hat\eta,x,\eps)= \hat J\vfi+ \eps^{-1}\hat\eta x+\eps S_1(\hat J, \vfi, \hat \eta, x)+\eps^2 S_2(\hat J, \vfi, \hat \eta, x),
 \end {equation}
 where $S_2$ is $2\pi$-periodic in $\vfi$ with average 0.
 We have:
\begin{equation}
      \label{var_trans_S2}
      \begin{cases}
      I=\hat J+\eps\dfrac{\p S_1(\hat J, \vfi, \hat \eta, x)}{\p \vfi}+\eps^2\dfrac{\p S_2(\hat J, \vfi, \hat \eta, x)}{\p \vfi}\, (1)\\
      \hat \psi=\vfi+\eps\dfrac{\p S_1(\hat J, \vfi, \hat \eta, x)}{\p \hat J}+\eps^2\dfrac{\p S_2(\hat J, \vfi, \hat \eta, x)}{\p \hat J}\, (2)\\
      y=\hat\eta+\eps^2\dfrac{\p S_1(\hat J, \vfi, \hat \eta, x)}{\p x}+\eps^3\dfrac{\p S_2(\hat J, \vfi, \hat \eta, x)}{\p x}\, (3)\\
      \hat \xi=x+\eps^2\dfrac{\p S_1(\hat J, \vfi, \hat \eta, x)}{\p \hat \eta}+\eps^3\dfrac{\p S_2(\hat J, \vfi, \hat \eta, x)}{\p \hat \eta}\, (4)
      \end{cases}.
\end{equation}
 The function $S_2$ should be chosen in such a way that  the Hamiltonian for the new variables has the form
\begin{equation}
      \label{H_transformed_S2}
     \hat{{\cal H}}(\hat J, \hat\psi, \hat\eta, \hat\xi, \eps)
     = H_0(\hat J,\hat\eta, \hat\xi)+ \eps {\cal H}_1(\hat J, \hat\eta, \hat\xi)
      + \eps^2 \hat{\cal H}_{2}(\hat J, \hat\eta, \hat\xi)+ \eps^3\hat{\cal H}_{3}(\hat J, \hat\psi, \hat\eta, \hat\xi, \eps).
\end{equation}
Then $H(I,\vfi,y,x,\eps)=\hat{\cal H}(\hat J, \hat\psi, \hat\eta, \hat\xi, \eps)$. Substitute (1), (3) of (\ref{var_trans_S2}) into $H(I,\vfi,y,x,\eps)$,  substitute (2), (4) of (\ref{var_trans_S2}) into $\hat{\cal H}(\hat J, \hat\psi, \hat\eta, \hat\xi, \eps)$, and equate terms of the same order in  expansions of these functions in $\eps$. We get
\begin{eqnarray}
   \label{H1H2_transformed_S2}
      {\cal{ H}}_1(\hat J, \hat\eta, x)&=& \dfrac{\p S_1}{\p \vfi}\dfrac{\p H_0(\hat J, \hat\eta, x)}{\p\hat J}+ H_1(\hat J, \vfi, \hat\eta, x,0),\\
      \hat{\cal H}_2(\hat J, \hat\eta, x)&=&\dfrac{\p S_2}{\p \vfi}\dfrac{\p H_0(\hat J, \hat\eta, x)}{\p\hat J}
      -\dfrac{\p S_1}{\p \eta}\dfrac{\p H_0(\hat J, \hat\eta, x)}{\p x}+\dfrac{\p S_1}{\p x}\dfrac{\p H_0(\hat J, \hat\eta, x)}{\p \hat \eta}
     \\
      &+&\dfrac 12\left(\dfrac{\p S_1}{\p \vfi}\right)^2\dfrac{\p^2 H_0(\hat J, \hat\eta, x)}{\p \hat J^2}+\dfrac{\p S_1}{\p \vfi}\dfrac{\p H_1(\hat J, \vfi, \hat\eta, x,0)}{\p \hat J}+\dfrac{\p H_1(\hat J, \vfi, \hat\eta, x,0)}{\p \eps}\nonumber.
\end{eqnarray}
The choice of $S_1, {\cal H}_1$ in Section \ref{s_main_formula}  ensures that the first equation in (\ref {H1H2_transformed_S2}) is satisfied.
Since  $S_1,S_2$ are $2\pi$-periodic functions of $\vfi$ with average 0, we integrate $\hat{\cal H}_2$ from 0 to $2\pi$ and get:
\begin{eqnarray*}
                             \label{H1H2S1S2_S2}
      \hat{\cal H}_2(\hat J, \hat\eta, x)&=&\dfrac{1}{2\pi}\int_{0}^{2\pi}\left(\dfrac 12\left(\dfrac{\p S_1}{\p \vfi}\right)^2\dfrac{\p^2 H_0}{\p \hat J^2}\right.
      +\left.\dfrac{\p S_1}{\p \vfi}\dfrac{\p H_1}{\p \hat J}+\dfrac{\p H_1}{\p \eps}\right) \, d\vfi,
      \end{eqnarray*}
      \vskip -0.5cm
      \begin{eqnarray}
      S_2(\hat J, \vfi, \hat \eta, x)&=&\frac{1}{\omega_0}\int_0^{\vfi} \left(\hat{\cal H}_2+\dfrac{\p S_1}{\p \eta}\dfrac{\p H_0}{\p x}-\dfrac{\p S_1}{\p x}\dfrac{\p H_0}{\p \hat \eta}\right.\\
      &-&\left.\dfrac 12\left(\dfrac{\p S_1}{\p \vfi}\right)^2\dfrac{\p^2 H_0}{\p \hat J^2}-\dfrac{\p S_1}{\p \vfi}\dfrac{\p H_1}{\p \hat J}-\dfrac{\p H_1}{\p \eps}\right)\, d\vfi_1 +s_2 (\hat J, \hat\eta, x)\nonumber,
\end{eqnarray}
where the function $s_2(\hat J, \hat\eta, x)$ is such that the average of $S_2$ is $0$.  One should put $\eps=0$ in right hand sides of (\ref{H1H2S1S2_S2}).

\begin{lemma}
\label{invert_1}
There exist positive constants $c_1, \ldots  , c_4$ such that the following is satisfied.\footnote{ The domain $D$ was introduced in Section \ref {s_assumptions}.}

(a) In the domain
\begin{equation}
\label{D_01}
(I,y, x,\vfi )\in  (D-c_1\sqrt{\eps})\times \Sm^1, \                                |\omega_0(I,y,x)|\ge c_2\sqrt {\eps}
    \end{equation}
 formulas  (\ref{transform_1})  and  (\ref {var_trans_S2}) determine transformations of variables $(I,\vfi,y,x)\mapsto (J,\psi,\eta, \xi)$
 and $(I,\vfi,y,x)\mapsto (\hat J,\hat \psi,\hat \eta, \hat \xi)$
 such that
 \begin{eqnarray}
\label{est_Iphiyx}
 I=J+O\Big(\frac{\eps}{\bar \omega_0}\Big), \ \vfi=\psi +O\Big(\frac{\eps}{\bar\omega_0^2}\Big), \ y=\eta+O\Big(\frac{\eps^2}{\bar\omega_0^2}\Big),\ x=\xi+O\Big(\frac{\eps^2}{\bar\omega_0^2}\Big),\\
 I=\hat J+O\Big(\frac{\eps}{\hat \omega_0}\Big), \ \vfi=\hat\psi +O\Big(\frac{\eps}{\hat \omega_0^2}\Big), \ y=\hat\eta+O\Big(\frac{\eps^2}{\hat \omega_0^2}\Big),\ x=\hat\xi+O\Big(\frac{\eps^2}{\hat\omega_0^2}\Big),\\
\label{J_Jhat}
  J=\hat J+O\Big(\frac{\eps^2}{\hat \omega^3_0}\Big), \ \psi=\hat\psi +O\Big(\frac{\eps^2}{\hat \omega_0^4}\Big), \ \eta=\hat\eta+O\Big(\frac{\eps^3}{\hat \omega_0^4}\Big),\ \xi=\hat\xi+O\Big(\frac{\eps^3}{\hat \omega_0^4}\Big).
  \end{eqnarray}
    Here  $\bar\omega_0=\omega_0(J,\eta,\xi), \ \hat \omega_0=\omega_0(\hat J,\hat \eta,\hat \xi)$.

 (b) The Hamiltonian in the new variables has forms (\ref{H_transformed_S1}),   (\ref{H_transformed_S2}) with the following estimates for  functions ${\cal H}_2, \hat{\cal H}_3$:
 \begin{eqnarray}
\label{etaxi_H_2}
 \frac{\p {\cal H}_2}{\p J}=O\Big(\frac{1}{\bar \omega_0^3}\Big),\ \frac{\p {\cal H}_2}{\p \psi}=O\Big(\frac{1}{\bar\omega_0^2}\Big), \ \frac{\p {\cal H}_2}{\p \eta}=O\Big(\frac{1}{\bar\omega_0^3}\Big),\  \frac{\p {\cal H}_2}{\p \xi}=O\Big(\frac{1}{\bar\omega_0^3}\Big),\\
\label{cal_H_3}
 \frac{\p\hat{\cal H}_3}{\p \hat J}=O\Big(\frac{1}{\hat\omega_0^5}\Big),\ \frac{\p\hat{\cal H}_3}{\p \hat \psi}=O\Big(\frac{1}{\hat\omega_0^4}\Big), \ \frac{\p\hat{\cal H}_3}{\p \hat \eta}=O\Big(\frac{1}{\hat\omega_0^5}\Big),\  \frac{\p\hat{\cal H}_3}{\p \hat\xi}=O\Big(\frac{1}{\hat\omega_0^5}\Big).
\end{eqnarray}
(c) Domains
\begin{equation}
\label{D_11}
(J,\psi,\eta, \xi )\in  (D-c_3\sqrt{\eps})\times \Sm^1, \                                |\omega_0(J,\eta, \xi )|\ge c_4\sqrt {\eps}
\end{equation}
and
\begin{equation}
\label{D_12}
(\hat J,\hat\psi,\hat\eta, \hat\xi )\in  (D-c_3\sqrt{\eps})\times \Sm^1, \    |\omega_0(\hat J,\hat\eta, \hat\xi )|\ge c_4\sqrt {\eps}
\end{equation}
belong to the image of domain (\ref{D_01}) under the considered transformation.

\medskip
(d)

\vskip-0.8cm
\begin{equation}
\begin{aligned}
&|\omega_0(I,y, x )|>\frac12|\omega_0(J,\eta, \xi )|, \ |\omega_0(I,y, x )|>\frac12|\omega_0(\hat J, \hat \eta, \hat \xi )|,\\
&  |\omega_0(J,\eta, \xi )|> \frac12|\omega_0(\hat J, \hat \eta, \hat \xi )|.
\end{aligned}
\end{equation}
\end{lemma}
This lemma is similar to standard assertions about the averaging method transformation (cf., e.g, \cite{Arn_mech}, Sec. 52C). We added an explicit  indication of orders of singularities at resonance (for  small $\omega_0$). We omit the proof.

\medskip

The phase point $(I(t), \vfi(t), y(t), x(t))$ belongs to the domain (\ref{D_01}) for $\tau$ close enough to $\tau_0$.   For such values of $\tau$ one can introduce  $( J(t), \psi(t), \eta(t),  \xi(t))$ and  $(\hat J(t),\hat\psi(t),\hat \eta(t), \hat \xi(t))$ via the transformations in Lemma \ref{invert_1}. \\
Denote $\omega_e(\tau)=\omega_0( I(t), y(t),  x(t)), \quad\bar\omega_0(\tau)=\omega_0( J(t), \eta(t),  \xi(t)), \\ \hat\omega_0(\tau)=\omega_0(\hat J(t),\hat \eta(t), \hat \xi(t))$.
For $\tau_1$  close enough  to $\tau_0$ on the time interval $\tau_0 \leq\tau\le\tau_1$, we have
\begin{equation}
\label{domain}
(I(t),y(t), x(t))\in  D-2c_1\sqrt{\eps}
\end{equation}
and
\begin{equation}
\label{condition_tau1}
 \omega_e(\tau)\ge\frac 12\omega_a(\tau), \ \bar\omega_0(\tau)\ge\frac 12\omega_a(\tau), \ \hat\omega_0(\tau)\ge\frac 12\omega_a(\tau), \ \omega_a(\tau)\ge 2c_2 \sqrt{\eps}.
\end{equation}
Denote $\tau_{1,sup}$ the supremum of values $\tau_1$ such that  relations (\ref{domain}), (\ref{condition_tau1}) are satisfied on the time interval $\tau_0\leq\tau\leq\tau_1$.
Consider motion on the time interval $\tau_0 \leq\tau\le\tau_{1,sup}$.  On this time interval the assumptions (\ref{D_01}) of Lemma  \ref{invert_1} are satisfied, and thus  estimates in this Lemma can be used.

\medskip
Estimates (\ref{line_3_6.1}),  (\ref{line_4_6.1}) of Lemma \ref{lemma_before} for this time interval follow from Lemma \ref{invert_1}
and the second condition in (\ref{condition_tau1}). Thus, we have
\begin{eqnarray}
\label{IJ_phipsi}
| I(t)-J(t)|=O\left(\frac{\eps}{\omega_a(\tau)}\right),\quad  |\vfi(t)-\psi(t)|=O\left(\frac{\eps}{\omega^2_a(\tau)}\right), \\
\label{yeta_xxi}
| y(t)-\eta(t)|=O\left(\frac{\eps^2}{\omega^2_a(\tau)}\right),\quad  |x(t)-\xi(t)|=O\left(\frac{\eps^2}{\omega^2_a(\tau)}\right).
\end{eqnarray}
\medskip
According to (\ref{cal_H_3}),
\begin{equation}
\dot {\hat J}= - \eps^3\dfrac{\p \hat{\cal H}_3}{\p \hat\psi}=O\left(\frac{\eps^3}{\hat\omega_0^4}\right)=O\left(\frac{\eps^3}{\omega_a^4}\right).
\end{equation}
Thus
\begin{equation}
\hat J(t)-\hat J_0=\int_{t_0}^t  O\left(\frac{\eps^3}{\omega_a^4(\eps \theta)}\right)d\theta=\int_{\tau_0}^{\tau}  O\left(\frac{\eps^2}{\omega_a^4(\vte)}\right)d\vte=     O  \left(\frac{\eps^2}{\omega_a^3(\tau)}\right).
\end{equation}
We use here that, according to condition (\ref{transvers}), solutions of the averaged system cross  the resonant surface transversally,  and thus ${d\omega_a(\tau)}/{d\tau}$ does not vanish when $\omega_a(\tau) =0$.
According to (\ref{J_Jhat}),
\begin{equation}
\label{JJhat_2}
\left| J(t)-\hat J(t)\right|=O\left(\frac{\eps^2}{\hat\omega_0^3}\right)=O\left(\frac{\eps^2}{\omega_a^3}\right).
\end{equation}
Then
\begin{equation}
\label{J_J_0}
\left| J(t)-J_0\right| \le \left| J(t)-\hat J(t)\right|+ \left| \hat J(t)-\hat J(t_0)\right| + \left| \hat J(t_0)-J_0\right| = O\left(\frac{\eps^2}{\omega_a^3}\right).
\end{equation}
Estimate $|\eta(t)-\eta_a(\tau)|,   |\xi(t)-\xi_a(\tau)|$. Equations for $\eta, \xi$ are determined by Hamiltonian (\ref{H_transformed_S1}):
 \begin{equation*}
 \begin{aligned}
 \dot \eta=-\eps\frac{\p {\cal H}_a(J, \eta,  \xi,\eps) }{\p \xi}-\eps^3\frac{\p {\cal H}_2(J, \psi, \eta,  \xi,\eps) }{\p \xi}  ,\\  \dot \xi=\eps\frac{\p {\cal H}_a(J, \eta,  \xi,\eps) }{\p \eta} +\eps^3\frac{\p {\cal H}_2(J, \psi, \eta,  \xi,\eps) }{\p \eta},
 \end{aligned}
  \end{equation*}
  where ${\cal H}_a=H_0+\eps {\cal H}_1$. Replace $J$ with $J_0$ in arguments of ${\cal H}_a$. This gives additional terms
   $O({\eps^3}/{ \omega_a^3})$ in the right hand sides of equations.  Using estimates for derivatives of $ {\cal H}_2$ in Lemma \ref{invert_1} we get that differential equations for $\eta,\xi$ coincide with differential equations for $\eta_a,\xi_a$ with an accuracy $O({\eps^3}/{ \omega_a^3})$. This implies estimates
  $$
  \eta(t)-\eta_a(\tau)=O\left(\frac{\eps^2}{ \omega_a^2}\right), \ \xi(t)-\xi_a(\tau)=O\left(\frac{\eps^2}{ \omega_a^2}\right).
  $$
  Together with already obtained estimates (\ref{IJ_phipsi}), (\ref{yeta_xxi}), (\ref{J_J_0}) this implies
  \begin{equation}
  \label{I_J}
   I(t)=J_0+O\left(\frac{\eps}{\omega_a}\right), \  y(t)=\eta_a(\eps t)+ O\left(\frac{\eps^2}{ \omega_a^2}\right),\  x(t)=\xi_a(\eps t)+ O\left(\frac{\eps^2}{ \omega_a^2}\right).
\end{equation}
Similarly, according to (\ref{J_Jhat}), (\ref{JJhat_2}) and (\ref{J_J_0}),
$$
 \hat J(t)=J_0+O\left(\frac{\eps^2}{\omega_a^3}\right), \  \hat \eta(t)=\eta_a(\eps t)+ O\left(\frac{\eps^3}{\omega_a^4}\right),\  \hat\xi(t)=\xi_a(\eps t)+ O\left(\frac{\eps^3}{\omega_a^4}\right).
$$

Then, we have
  $$
  \dot \psi=\frac{\p {\cal H}_a(J, \eta,  \xi,\eps) }{\p J}+\eps^2\frac{\p {\cal H}_2(J, \psi, \eta,  \xi,\eps) }{\p J}.
  $$
   Replace $J, \eta,\xi $ with $J_0,  \eta_a,\xi_a$ in arguments of ${\cal H}_a$. This gives additional terms
   $O({\eps^2}/{\omega_a^3})$ in the right hand sides of formula for $\dot \psi$.  Using estimates for derivatives of $ {\cal H}_2$ in Lemma \ref{invert_1} we get,
   $$
   \dot \psi= \dot\psi_a+O\left(\frac{\eps^2}{\omega_a^3}\right).
   $$
   This implies
   $$
   \psi(t)=\psi_a(\tau)+O\left(\frac{\eps}{\omega_a^2}\right).
   $$

   \medskip
    On the time interval $\tau_0 \leq\tau\le\tau_{1,sup}$ we have
    \begin{equation}
\label{domain+}
(I(t),y(t), x(t))\in  D-3c_1\sqrt{\eps}.
\end{equation}
This is because the solution of the averaged system with the initial condition $(I_0,  y_0, x_0)$ stays at a distance of order 1  from the boundary of $D$ (due to assumption D in Section \ref{s_assumptions}),  and the component $(I(t), y(t),x(t))$  of the solution of system (\ref {e_system}) with the initial condition  $(I_0, \vfi_0, y_0, x_0)$ is at a distance $O(\sqrt{\eps})$ from this solution of the averaged system (as it follows from estimates (\ref{I_J})). Therefore, condition (\ref{domain}) is satisfied on this time interval with a margin and is still valid for some time after  $\tau_{1,sup}$.

   \medskip
      On the time interval $\tau_0 \le \tau\le\tau_{1,sup}$ we have
   \begin{eqnarray*}
   \omega_e(\tau)=\omega_0( I(t),y(t),  x(t))= \omega_0( J_0,\eta_a(\tau),  \xi_a(\tau)) + O\left(\frac{\eps}{\omega_a}\right)=\omega_a(\tau)+ O\left(\frac{\eps}{\omega_a(\tau)}\right),\\
\bar\omega_0(\tau)=\omega_0( J(t),\eta(t),  \xi(t))= \omega_0( J_0,\eta_a(\tau),  \xi_a(\tau)) + O\left(\frac{\eps^2}{\omega_a^3}\right)=\omega_a(\tau)+ O\left(\frac{\eps}{\omega_a(\tau)}\right),\\
\hat\omega_0(\tau)=\omega_0(\hat J(t),\hat \eta(t), \hat \xi(t))= \omega_0( J_0,\eta_a(\tau),  \xi_a(\tau))
+ O\left(\frac{\eps^2}{\omega_a^3}\right)=\omega_a(\tau)+O\left(\frac{\eps}{\omega_a(\tau)}\right).
\end{eqnarray*}
Therefore there exists a positive constant $d_1$ such that
$$
\omega_e(\tau)>\omega_a(\tau)-  d_1\frac{\eps}{\omega_a(\tau)}, \ \bar\omega_0(\tau)>\omega_a(\tau)-  d_1\frac{\eps}{\omega_a(\tau)}, \  \hat\omega_0(\tau)>\omega_a(\tau)-  d_1\frac{\eps}{\omega_a(\tau)}.
$$
Choose constant $C_1$ in the formulation of  Lemma \ref{lemma_before} such that $C_1>2c_2$ and $C_1> 2\sqrt{d_1}$. Then for $ \omega_a(\tau)\ge C_1\sqrt{\eps}$ we have
$d_1{\eps}/{\omega_a(\tau)}<{\omega_a(\tau)}/{4}$ and, therefore,
\begin{equation}
\label{for_3/4}
\begin{aligned}
\omega_e(\tau)>\omega_a(\tau)-\frac{\omega_a(\tau)}{4}=\frac{3\omega_a(\tau)}{4},\\
\bar\omega_0(\tau)>\omega_a(\tau)-\frac{\omega_a(\tau)}{4}=\frac{3\omega_a(\tau)}{4}, \\
\hat \omega_0(\tau)>\omega_a(\tau)-\frac{\omega_a(\tau)}{4}=\frac{3\omega_a(\tau)}{4}.
\end{aligned}
\end{equation}
Denote $\tau_{1,a}$ the moment of slow time such that  $ \omega_a(\tau_{1,a})= C_1\sqrt{\eps}$.  Estimates (\ref{for_3/4}) imply that $\tau_{1,a}< \tau_{1,sup}$.  Therefore,   assumptions  (\ref{condition_tau1}) are satisfied, and obtained estimates of Lemma \ref{lemma_before} are valid while $ \omega_a(\tau)\ge C_1\sqrt{\eps}$. This completes the proof of Lemma \ref{lemma_before}.

\subsection{Proof of Lemma \ref{lemma_after}}
\label{proof_lemma_after}
Statement (a) of  Lemma~\ref{lemma_after} is a direct corollary of Lemma  \ref{lemma_before} and conditions C and D from Section~\ref{s_assumptions}.

For the statements (b) and (c) we need the following statement.
\begin{lemma} \label{l:pass-res}
  For any $r>1$ there exists an exceptional set $\mathcal V_r$ with measure $O(\varepsilon^r)$ such that the following holds. All solutions of the perturbed system with initial data outside $\mathcal V_r$ reach $\omega_0 = 0$ at some moment $t_e > t_N$.
  If $F$ has critical points,
  \begin{equation} \label{e:t-N-e}
    t_e - t_N = O\Big(\frac{\ln \nu}{\sqrt{\varepsilon}}\Big)
  \end{equation}
  with $\nu \ge \varepsilon$ depending on the initial data as described in the statement of Lemma~\ref{l:pass-res}). If $F$ has no critical points,
  \begin{equation}
    \varepsilon(t_e - t_N) = O(\sqrt{\varepsilon}).
  \end{equation}
\end{lemma}
Let us explain how Lemma~\ref{l:pass-res} implies statements (b) and (c) before proving it. We will do this for statement (b), for (c) it can be done similarly.

\begin{enumerate}
  \item As $\varepsilon t_N = \tau_* + O(\sqrt{\varepsilon})$, \eqref{e:t-N-e} implies $\varepsilon t_e = \tau_* + O(\sqrt{\varepsilon} \ln \nu)$.
  \item According to Lemma~\ref{lemma_before},
  \begin{equation} \label{e:tN}
    |I(t_N) -I_0|= O(\sqrt{\eps}),\;  |y(t_N) -\eta_a(\eps t_{N})|= O({\eps}), \; |x(t_N) -\xi_a(\eps t_{N})|= O({\eps}).
  \end{equation}
  \item As $\dot I = O(\varepsilon)$, \eqref{e:t-N-e} and~\eqref{e:tN} imply $I(t) = I_0 + O(\sqrt{\varepsilon} \ln \nu)$ for all $t \in [t_N, t_e]$.
  \item We have
  \begin{align}
  \begin{split}
    \frac{dy}{dt} - \frac{d}{dt} \eta_a(\varepsilon t_N) &=
    -\varepsilon \Big( \frac{\partial H_0}{\partial x}\Big|_{I, y, x} - \frac{\partial H_0}{\partial x}\Big|_{I_0, \eta_a, \xi_a} \Big) + O(\varepsilon^2) =\\
    &= O(\varepsilon \sqrt{\varepsilon} \ln \nu).
  \end{split}
  \end{align}
\end{enumerate}
We have used that $I-I_0 = O(\sqrt{\varepsilon} \ln \nu)$.
Together with~\eqref{e:t-N-e} and~\eqref{e:tN} this estimate implies $y(t) = \eta_a(\varepsilon t) + O(\varepsilon \ln^2 \nu)$. We get $x(t) = \xi_a(\varepsilon t) + O(\varepsilon \ln^2 \nu)$ in the same way.
We have proved Lemma~\ref{lemma_after} using Lemma~\ref{l:pass-res}.

Let us now prove Lemma~\ref{l:pass-res}.
We will use the rescaled action $P=(I-a(y,  x))/\sqrt{\eps}$ and the rescaled time $\theta = \sqrt{\eps} t$ introduced in Section~\ref{s_pendulum} as well as the function ${\cal E}$ given by formula~\eqref{def_E1}.
Consider the motion given by the perturbed system~\eqref{a_system} in the chart $(P, \vfi, \tilde y, \tilde x)$ introduced in Section~\ref{s_pendulum}.
This motion writes as (cf. Section~\ref{s_pendulum}).
\begin{equation*}
  \dot P = - \frac{\partial F}{\partial \vfi}(\tilde y, \tilde x) + O(\sqrt{\varepsilon}),
  \quad
  \dot \vfi = P + O(\sqrt{\varepsilon}),
  \quad
  \dot {\tilde y} =  O(\sqrt{\varepsilon}), \quad \dot {\tilde x} =  O(\sqrt{\varepsilon}).
\end{equation*}
The principal part is the Hamiltonian system with the Hamiltonian $\mathcal E(P, \vfi)$:
\begin{equation} \label{e:particle}
  \dot P = - \frac{\partial F}{\partial \vfi}(\tilde y, \tilde x),
  \qquad
  \dot \vfi = P,
  \qquad
  \dot {\tilde y} = 0, \qquad \dot {\tilde x} = 0.
\end{equation}
Let us start with the case when the function $F(\vfi, y_*,x_*)$ does not have critical points. Then the projection of the phase point $(I(t), \vfi(t), y(t), x(t))$ onto the plane $P, \vfi$ moves approximately along a level line of the Hamiltonian ${\cal E}(P, \vfi,   y_*,x_*)$ in the phase portrait~\ref{portraits}a and arrives at the resonance $P=0$ at time $t_e= t_{N}+O(1/\sqrt{\eps}\,)$. Indeed, we have $\frac{dP}{d\theta} \approx -\frac{\partial F}{\partial \varphi} < 0$ and $P(t_N) = O(1)$. Clearly, estimates (\ref{est_br_2}) are satisfied. This proves Lemma~\ref{l:pass-res} for the case without critical points.

Let us now consider the case  when $F(\vfi, y_*,x_*)$ has critical points.
Critical points of $F$ correspond to fixed points of~\eqref{e:particle} restricted to the $(P,\varphi)$-plane (with $P=0$), and there are no other fixed points. Local maxima of $F$ correspond to saddles, and in the full phase space of~\eqref{e:particle} these saddles become codimension two normally hyperbolic invariant manifolds having a codimension one stable and unstable manifolds that correspond to the separatrices of the saddles.
By Fenichel's results these normally hyperbolic manifolds survive a small perturbation, and near the continuation of each of these manifolds there exists a linearizing chart~\cite{Fen}.
More precisely, let $\vfi_i$ be a local maxima of $F(\vfi, x_*, y_*)$, then there are new fast variables $(A, B)$ defined in some neighbourhood of the point $(\vfi_i, 0, x_*, y_*)$ such that the dynamics of $A$ and $B$ given by the perturbed system~\eqref{a_system} writes as
\begin{align}
\begin{split} \label{e:Fenichel}
  \dot A &= -\lambda(\tilde y, \tilde x) A \Big(1+O(|A|+|B|+|\sqrt{\eps}|)\Big), \\
  \dot B &= \mu(\tilde y, \tilde x) B  \Big(1+O(|A|+|B|+|\sqrt{\eps}|)\Big),
\end{split}
\end{align}
where $\lambda, \mu > 0$.
We see that the stable and the unstable manifolds write as $A=0$ and $B=0$, respectively.
Trajectories passing far from all normally hyperbolic manifolds described above reach $P=0$ after time $\theta$ of order $O(1)$, so we have $\varepsilon t_e = \varepsilon t_N + O(\sqrt{\varepsilon})$ for such trajectories.
However, the motion slows near these normally hyperbolic manifolds, so extra work is required to prove Lemma~\ref{l:pass-res}.

Denote by $\{S_i\}$ the normally hyperbolic manifolds of~\eqref{a_system}. Recall that we are assuming that the values of $F$ at its critical points are different from each other (Condition~E). Then the values of $\mathcal E$ are also different. Together with the fact that $\mathcal E$ is almost conserved, this means that each trajectory can come close to at most one of the manifolds $S_i$.
This condition also implies that stable and unstable manifolds of different manifolds $S_i$ are disjoint.

Cover each manifold $S_i$ by a small neighborhood $U_i$ of the form $|A| < c, |B| < c$, where $c > 0$ is so small that the values of $\mathcal E$ are disjoint in different $U_i$. Then solutions of~\eqref{a_system} can intersect at most one neighborhood $U_i$. Consider a solution visiting $U_i$.
Let $B_{in}$ be the value of $B$ at the time of the first entry into $U_i$. The continuous time $\theta$ spent inside $U_i$ during this entry is $O(|\ln B_{in}|)$, as can be checked using~\eqref{e:Fenichel}.
So, for any $r>1$ there exists $c_r > 0$ such that time spent in this neighborhood is greater than $c_r \ln \varepsilon$ only when $|B_{in}| < \varepsilon^r$. This inequality on $B_{in}$ means that the $(P, \varphi)$-distance between this solution and the stable manifold of $S_i$ is $\lesssim \varepsilon^r$ when $P \sim 1$. Set as $\mathcal V_r$ the set of initial data corresponding to such solutions. The measure of $\mathcal V_r$ is $O(\varepsilon^r)$. Solutions with initial data outside $\mathcal V_r$ spend time $\theta$ at most $O(\ln \varepsilon)$ in $U_i$, and then leave it while being near a separatrix loop of~\eqref{e:particle} (cf. Figure~\ref{portraits}b) and reach $P=0$. This gives the estimate $\varepsilon t_e = \varepsilon t_N + O(\sqrt{\varepsilon} \ln \varepsilon)$ for all initial data outside $\mathcal V_r$.

To get the estimate $\varepsilon t_e = \varepsilon t_N + O(\sqrt{\varepsilon} \ln \nu_1)$, we may only consider the case $\nu_1 \gg \sqrt{\varepsilon} |\ln \varepsilon|$ (otherwise, $\ln \nu \sim \ln \eps$).
As $\mathcal E$ is a first integral of~\eqref{e:particle}, we get $\frac{d \mathcal E}{d\theta} = O(\sqrt{\varepsilon})$ and $\mathcal E(t) = \mathcal E(t_e) + O(\sqrt{\varepsilon} \ln \varepsilon)$ for $t \in [t_N, t_e]$. By the definition of $\nu_1$ we have $\mathcal E(t_e) = \mathcal E_{crit} + O(\nu_1^2)$, where $\mathcal E_{crit}$ is the value of $F$ at the critical point corresponding to $S_i$. This means that at the moment of the entry in $U_i$ we have $|B_{in}| \sim \nu_1^2$, and the time spent in $U_i$ is $O(\ln B_{in}) = O(\ln \nu_1)$.

\section {Proofs of lemmas about estimates of integrals}
\label{est_int}

{\bf Proof of Lemma \ref {last_integral}.}
We have $t_e=t_N+O\left(\frac{1}{\sqrt{\eps}}\ln\nu\right)$.  All terms   in the integrand but the term proportional to $(b(y,x)-b_{*,a})$ are $O(\eps)$.
Their contribution in the integral is $O(\sqrt{\eps}\ln \nu)$. \footnote{Contributions of terms proportional to $(I-a(y,x))^2$ or $\eps(I-a(y,x))$ is  $O(\sqrt{\eps})$.}

 Dynamics for
$t_N\le t\le t_e$ is described by the Hamiltonian (see Section \ref{s_pendulum})
\begin{equation}
\label{def_E2}
{\cal E}_{\tilde y,\tilde  x,\eps}(P, \vfi)=\frac 12\alpha(\tilde y, \tilde x)P^2+F(\vfi,  \tilde y,  \tilde x) +O(\sqrt{\eps})P+O(\eps).
\end{equation}
Here $P$ and  $\vfi $ are conjugate canonical variables, the time is
$\theta =\sqrt {\eps} t$, variables $ \tilde y, \tilde x$ change with the speed $O(\sqrt{\eps})$ in this time.

Consider some small neighborhoods of the saddle points on the Figure~\ref{portraits}b, denote them by $U_i(\tilde x, \tilde y) \subset \mathbb R^2_{P, \varphi}$, where $i$ enumerates these saddles. The time $\theta$ spent outside $\cup_i U_i$ is $O(1)$, this corresponds to a change of $\varphi$ of order $O(1)$. The motion slows down near the saddle points.
Hovever, we have
\begin{equation}
  \frac{d\tilde x}{dt} = O(\varepsilon),
  \qquad
  \frac{d\tilde y}{dt} = O(\varepsilon),
  \qquad
  t_e - t_N = O\Big(\frac 1 {\sqrt \varepsilon} \ln \nu\Big),
\end{equation}
so the total change of $\tilde x$ and $\tilde y$ between $t_N$ and $t_e$ is $O(\sqrt{\varepsilon} \ln \nu) = O(1)$, and the saddle points also have moved by $O(\sqrt{\varepsilon} \ln \nu) = O(1)$.
So, the total change of $\varphi$ is $O(1)$.

Consider the term $\alpha(y,x)(I-a(y,x)) (b(y,x)-b_{*,a})$ in the integrand  in  Lemma \ref {last_integral}. We  have $ (b(y,x)-b_{*,a})= O(\sqrt{\eps}\ln \nu)$, and
$\alpha(y,x)(I-a(y,x))=\dot \vfi +O(\eps)$.    Change of $\vfi$ in the considered time interval is $O(1)$. Thus the contribution of the term $\alpha(y,x)(I-a(y,x)) (b(y,x)-b_{*,a})$  in the integral is $O(\sqrt{\eps}\ln \nu)$. This implies the result of the Lemma.

\bigskip
{ \bf Proof of Lemma \ref{int_with S1}.}
We use integration by parts:
 \begin{equation}
\begin{aligned}
&\eps\int_{0}^{t_N}\left(J_0-a(\eta_a,\xi_a)\right)\frac {\p S_1}{\p\vfi}L_1dt\\
=&\eps\int_{0}^{t_N}\left(     \frac {d S_1}{d t} -\eps (\frac {\p S_1}{\p y}\eta_a'+\frac {\p S_1}{\p x}\xi_a')
\right)\frac{1}{\dot \psi}\left(J_0-a(\eta_a,\xi_a)\right)L_1dt\\
=&\eps S_1 \frac{1}{\dot \psi}\left(J_0-a(\eta_a,\xi_a)\right)L_1(J_0,\eta_a,\xi_a) \vert_0^{t_N}\\
-&\eps^2\int_{0}^{t_N}S_1\frac{d}{d\tau}\left(\frac{1}{\dot \psi}\left(J_0-a(\eta_a,\xi_a)\right)L_1(J_0,\eta_a,\xi_a) \right)dt\\
-&\eps^2\int_{0}^{t_N}\left(\frac {\p S_1}{\p y}\eta_a'+\frac {\p S_1}{\p x}\xi_a'\right)
\frac{1}{\dot \psi}\left(J_0-a(\eta_a,\xi_a)\right)L_1(J_0,\eta_a,\xi_a)dt\\
=&O(\sqrt{\eps})+\int_{0}^{\tau_N}O\left(\frac{\eps}{\omega_a^2(\tau)}\right)d\tau=O(\sqrt{\eps}).
\end{aligned}
\end{equation}

\bigskip

{\bf Proof of Lemma \ref{lemma_T0}.}
We have
\begin{equation}
\begin{aligned}
\label{expan_lemma1n}
L=&(I-a(y,x))^2L_1(I,y,x)\\
=&\left(J_0-a(\eta_a,\xi_a)+\eps\frac {\p S_1}{\p\vfi}\left(J_0,\psi,\eta_a,\xi_a\right)
\right)^2L_1(J_0+\eps\frac {\p S_1}{\p\vfi},\eta_a,\xi_a)\\
+&O\left(\frac{\eps^2}{\omega^2_a(\tau)}\right)\\
=&(J_0-a(\eta_a,\xi_a))^2L_1(J_0,\eta_a,\xi_a)\\
+&2\eps\left(J_0-a(\eta_a,\xi_a)\right)\frac {\p S_1}{\p\vfi}\,L_1(J_0,\eta_a,\xi_a)+\eps(J_0-a(\eta_a,\xi_a))^2\frac {\p L_1}{\p I}\frac {\p S_1}{\p\vfi}\\
+& O\left(\frac{\eps^2}{\omega_a^2(\tau)}\right).
\end{aligned}
\end{equation}
We have
\begin{equation}
\begin{aligned}
\int_{0}^{t_N}(J_0-a(\eta_a,\xi_a))^2L_1(J_0,\eta_a,\xi_a)dt
&=\frac{1}{\eps}\int_{0}^{\tau_N}(J_0-a(\eta_a,\xi_a))^2L_1(J_0,\eta_a,\xi_a)d\tau\\
=\frac{1}{\eps}\int_{0}^{\tau_{*,a}}(J_0-a(\eta_a,\xi_a))^2L_1(J_0,\eta_a,\xi_a)d\tau &+O(\sqrt{\eps}),
\end{aligned}
\end{equation}
where $\tau_N=\eps t_N$.
Note that   $\tau_{*,a}- \tau_N= O(\sqrt{\eps})$, and      the integrand is $O(\eps)$ for $\tau_N \le\tau\le\tau_{*,a}$.

The integral of terms containing $(\p S_1/ \p \vfi)$ is $O(\sqrt{\eps})$ according to Lemma \ref{int_with S1}.
\begin{equation}
\begin{aligned}
&\int_{0}^{t_N} O\left(\frac{\eps^2}{\omega_a^2(\tau)}\right)dt=O({\sqrt \eps}).
\end{aligned}
\end{equation}
Combining all the estimates we get the result of the Lemma.

\bigskip

{\bf Proof of Lemma \ref{lemma_T1}.}
Expand $L$ using estimates (\ref{Part8_con_1}):
\begin{equation}
\begin{aligned}
\label{expan_lemmaT1}
&{L}(I,y,x,\eps)= L(J_0,\eta_a,\xi_a,0)+\eps\frac {\p S_1}{\p\vfi}\left(J_0,\psi,\eta_a,\xi_a\right)\frac {\p L}{\p I}\left(J_0,\eta_a,\xi_a,0\right)\\
+&O\left(\frac{\eps^2}{\omega_a^3(\tau)}\right)+O(\eps).
\end{aligned}
\end{equation}
We have
\begin{equation*}
\begin{aligned}
\eps\int_{0}^{t_N}L(J_0,\eta_a,\xi_a,0)dt
=\int_{0}^{\tau_{*,a}}L(J_0,\eta_a,\xi_a,0)d\tau &+O(\sqrt{\eps}).
\end{aligned}
\end{equation*}
Integration by parts gives the estimate
\begin{equation*}
\begin{aligned}
\eps^2\int_{0}^{t_N}\frac {\p S_1}{\p\vfi}\left(J_0,\psi,\eta_a,\xi_a\right)\frac {\p L}{\p I}\left(J_0,\eta_a,\xi_a,0\right)dt=O({\eps}).
\end{aligned}
\end{equation*}
Also
\begin{equation}
\begin{aligned}
&\int_{0}^{t_N} \left(O\left(\frac{\eps^3}{\omega_a^3(\tau)}\right) +O(\eps^2)\right)dt=O( \eps).
\end{aligned}
\end{equation}
Combining all the estimates we get the result of the Lemma.

\bigskip

{\bf Proof of Lemma \ref{lemma_T2}.}
Expand $L$ using estimates (\ref{Part8_con_1}):
\begin{equation}
\begin{aligned}
\label{expan_lemmaT2}
&{L}(I,\vfi,y,x,\eps)= L(J_0,\psi,\eta_a,\xi_a,0)+\eps\frac {\p S_1}{\p\vfi}\left(J_0,\psi,\eta_a,\xi_a\right)\frac {\p L}{\p I}\left(J_0,\psi,\eta_a,\xi_a,0\right)\\
-&\eps\frac {\p S_1}{\p I}\left(J_0,\psi,\eta_a,\xi_a\right)\frac {\p L}{\p \vfi}\left(J_0,\psi,\eta_a,\xi_a,0\right)+O\left(\frac{\eps^2}{\omega_a^4(\tau)}\right)+O(\eps)\\
=&L(J_0,\psi,\eta_a,\xi_a,0)+O\left(\frac{\eps}{\omega_a^2(\tau)}\right)
.
\end{aligned}
\end{equation}
 Denote $k(I,\vfi,y,x)=\int_0^{\vfi} {L}(I,\beta,y,x,0) d\beta$. Note that $k(I,\vfi,y,x,\eps)$ is also a periodic function of $\vfi$. Then
 \begin{equation}
\begin{aligned}
\int_{0}^{t_N}\eps L(J_0,\psi,\eta_a,\xi_a,0)dt=\int_{0}^{t_N}\eps\frac{\p k}{\p \vfi}(J_0,\psi,\eta_a,\xi_a,0)dt.
\end{aligned}
\end{equation}
Integration by parts gives the estimate  $O(\sqrt\eps)$ for this integral.

Also
\begin{equation}
\begin{aligned}
&\int_{0}^{t_N} O\left(\frac{\eps^2}{\omega_a^2(\tau)}\right)dt=O(\sqrt \eps).
\end{aligned}
\end{equation}
Combining all the estimates we get the result of the Lemma.

\bigskip

{\bf Proof of Lemma \ref {lemma_T3}.}
Expand $L$ using estimates (\ref{Part8_con_1}):
\begin{equation}
\begin{aligned}
\label{expan_lemmaT3}
&{L}(I,y,x)= L(J_0,\eta_a,\xi_a)+\eps\frac {\p S_1}{\p \vfi}\left(J_0,\psi,\eta_a,\xi_a\right)
\alpha (\eta_a,\xi_a)(b(\eta_a,\xi_a)-b_{*,a})\\
+&O\left(\frac{\eps^2}{\omega_a^2(\tau)}\right).
\end{aligned}
\end{equation}
We have
\begin{equation}
\begin{aligned}
\int^{t_N}_0 L(J_0,\eta_a,\xi_a) dt=\frac{1}{\eps}\int_{0}^{\tau{*,a}}L(J_0,\eta_a,\xi_a)dt+
O(\sqrt\eps).
\end{aligned}
\end{equation}
The integral of the term containing $(\p S_1/ \p \vfi)$ is $O(\sqrt{\eps})$ as the integration by parts shows. Also
\begin{equation}
\begin{aligned}
&\int_{0}^{t_N} O\left(\frac{\eps^2}{\omega_a^2(\tau)}\right)dt=O({\sqrt \eps}).
\end{aligned}
\end{equation}
Combining all the estimates we get the result of the Lemma.

\bigskip

{ \bf Proof of Lemma \ref{lemma_T4}.}
Expand $H_1$  and $K$ using estimates (\ref{Part8_con_1}):
\begin{equation}
\begin{aligned}
\label{expan1_lemmaT4}
&\frac{\p H_1}{\p \vfi}(I,\vfi,y,x,\eps)=\frac{\p H_1}{\p \vfi}(J,\psi,\eta_a,\xi_a,\eps)
+\eps\frac {\p S_1}{\p \vfi}\left(J_0,\psi,\eta_a,\xi_a\right)\frac{\p^2 H_1}{\p \vfi\p I}(J,\psi,\eta_a,\xi_a,\eps)\\
-&\eps\frac {\p S_1}{\p I}\left(J_0,\psi,\eta_a,\xi_a\right)\frac{\p^2 H_1}{\p \vfi^2}(J,\psi,\eta_a,\xi_a,\eps)+O\left(\frac{\eps^2}{\omega_a^4(\tau)}\right),
\end{aligned}
\end{equation}
\begin{equation}
\begin{aligned}
\label{expan2_lemmaT4}
&{K}(I,y,x)= K(J_0,\eta_a,\xi_a)+\eps\frac {\p S_1}{\p \vfi}\left(J_0,\psi,\eta_a,\xi_a\right)
\frac {\p K}{\p I}(J_0,\psi, \eta_a,\xi_a)+O\left(\frac{\eps^2}{\omega_a^2(\tau)}\right).
\end{aligned}
\end{equation}
Then
\begin{equation}
\begin{aligned}
\label{expan3_lemmaT4}
& \frac{\p H_1}{\p \vfi}(I,\vfi,y,x,\eps) K(I,y,x)=\frac{\p H_1}{\p \vfi}(J,\psi,\eta_a,\xi_a,\eps)
 K(J_0,\eta_a,\xi_a)\\
 +&\eps\left(\frac {\p S_1}{\p \vfi}\frac{\p^2 H_1}{\p \vfi\p I}K- \frac {\p S_1}{\p I}\frac{\p^2 H_1}{\p \vfi^2}K
 +\frac{\p H_1}{\p \vfi}\frac {\p S_1}{\p \vfi}
\frac {\p K}{\p I} \right)
+O\left(\frac{\eps^2}{\omega_a^2(\tau)}\right).
\end{aligned}
\end{equation}
Here it is used that $K=O(\omega_a^2)$.

Denote ${\cal H}_{1,\eps}(I, y,x,\eps)$ the average of  $H_1(I,\vfi,y,x,\eps)$ over $\vfi$.
Denote
\[
  \tilde H_{1,\eps}(I,\vfi,y,x,\eps)=H_1(I,\vfi,y,x,\eps)-{\cal H}_{1,\eps}(I, y,x,\eps).
\]
We have
\begin{equation}
\begin{aligned}
\label{int1_lemmaT4}
&\int^{t_N}_0 \frac{\p H_1}{\p \vfi}(J,\psi,\eta_a,\xi_a,\eps)
 K(J_0,\eta_a,\xi_a)dt=\int^{t_N}_0 \frac{\p \tilde H_{1,\eps}}{\p \vfi}K dt\\
 =&\int^{t_N}_0\frac{1}{\dot \psi}\left(\frac{d\tilde H_{1,\eps}}{dt} -   \eps \left(\frac{\p \tilde H_{1,\eps}}{\p y}\eta_a'+ \frac{\p \tilde H_{1,\eps}}{\p x}\xi_a'\right)    \right)Kdt\\
 =&\tilde H_{1,\eps}\frac{K}{\dot \psi}|_0^{t_N}
 -\eps \int^{t_N}_0\tilde H_{1,\eps} \left(\frac{K}{\dot \psi}\right)'dt
 -\eps \int^{t_N}_0\left(\frac{\p \tilde H_{1,\eps}}{\p y}\eta_a'+ \frac{\p \tilde H_{1,\eps}}{\p x}\xi_a'    \right)\frac{K}{\dot \psi}dt.
\end{aligned}
\end{equation}
We have
\begin{equation}
\begin{aligned}
\tilde H_{1,\eps}\frac{K}{\dot \psi}|_0^{t_N}=
-\frac{\tilde H_1(I_0, \vfi_0,y_0,x_0)}{\omega_0(I_0,y_0,x_0)}K(I_0,y_0,x_0)+O({\sqrt \eps}).
\end{aligned}
\end{equation}
Integration by parts in (\ref{int1_lemmaT4}) shows that terms proportional to $\eps$  are
$O({\sqrt \eps})$ (cf. Lemma \ref{lemma_T2}, note that average of $\tilde H_{1,\eps}$ over $\psi$ equals 0).

Now we should estimate the integral from 0 to $t_N$ of terms with $S_1$ in the second line of (\ref{expan3_lemmaT4}). According to Lemma \ref{lemma_T2}, this integral can be replaced by the integral of averaged of over $\psi$ integrand with the accuracy $O({\sqrt \eps})$. And then this integral can be replaced by the integral from 0 to $\tau_{*,a}$ with the same accuracy. Denote averaging over $\psi$ by the overline. Thus we have
\begin{equation}
\begin{aligned}
\label{long_averaged}
&\eps\int^{t_N}_0 \left(\frac {\p S_1}{\p \vfi}\frac{\p^2 H_1}{\p \vfi\p I}K- \frac {\p S_1}{\p I}\frac{\p^2 H_1}{\p \vfi^2}K +\frac{\p H_1}{\p \vfi}\frac {\p S_1}{\p \vfi}
\frac {\p K}{\p I}  \right)dt\\
=&\int^{\tau_{*,a}}_0 \overline{\left(\frac {\p S_1}{\p \vfi}\frac{\p^2 H_1}{\p \vfi\p I}K- \frac {\p S_1}{\p I}\frac{\p^2 H_1}{\p \vfi^2}K +\frac{\p H_1}{\p \vfi}\frac {\p S_1}{\p \vfi}
\frac {\p K}{\p I}  \right)}d\tau + O({\sqrt \eps}).
\end{aligned}
\end{equation}
From now on we assume that  $\eps$ is replaced by 0 in arguments of  $H_1$. This does not change the accuracy.
 We can replace $H_1$ with $\tilde H_1$ in (\ref{long_averaged}). Recall that
$\p S_1/\p \vfi=- \tilde H_1/\omega_0$.
Then
\begin{equation}
\begin{aligned}
\overline{\frac{\p H_1}{\p \vfi}\frac {\p S_1}{\p \vfi}\frac {\p K}{\p I}}= -\frac{1}{\omega}\frac {\p K}{\p I}\overline{\frac{\p \tilde H_1}{\p \vfi}\tilde H_1}= -\frac{1}{2\omega}\frac {\p K}{\p I}\overline{\frac{\p {\tilde H_1}^2}{\p \vfi}}=0.
\end{aligned}
\end{equation}
Also
\begin{equation}
\begin{aligned}
 &\overline{\frac {\p S_1}{\p \vfi}\frac{\p^2 \tilde H_1}{\p \vfi\p I}- \frac {\p S_1}{\p I}\frac{\p^2 \tilde H_1}{\p \vfi^2} }
 = \overline{\frac {\p S_1}{\p \vfi}\frac{\p^2 \tilde H_1}{\p \vfi\p I}
 -\frac {\p S_1}{\p I}   \frac{\p }{\p \vfi}  \left(\frac{\p\tilde H_1}{\p \vfi} \right)}  \\
 =&  \overline{\frac {\p S_1}{\p \vfi}\frac{\p^2 \tilde H_1}{\p \vfi\p I}
 -\frac {\p S_1}{\p I}   \frac{\p }{\p \vfi}  \left(\frac{\p\tilde H_1}{\p \vfi} \right)
 - \frac{\p }{\p \vfi}\left(\frac {\p S_1}{\p I}\right)   \frac{\p\tilde H_1}{\p \vfi} +\frac{\p }{\p \vfi}\left(\frac {\p S_1}{\p I}\right)   \frac{\p\tilde H_1}{\p \vfi}}\\
 =& \overline{\frac {\p S_1}{\p \vfi}\frac{\p^2 \tilde H_1}{\p \vfi\p I}
 - \frac{\p }{\p \vfi} \left(\frac {\p S_1}{\p I}\frac{\p \tilde H_1}{\p \vfi}  \right)
 +\frac{\p }{\p \vfi}\left(\frac {\p S_1}{\p I}\right)   \frac{\p\tilde H_1}{\p \vfi}}\\
 =& \overline{\frac {\p S_1}{\p \vfi}\frac{\p^2 \tilde H_1}{\p \vfi\p I}
 +\frac{\p }{\p \vfi}\left(\frac {\p S_1}{\p I}\right)   \frac{\p\tilde H_1}{\p \vfi}}
 =\overline{\frac{\p }{\p I}\left(\frac {\p S_1}{\p \vfi}\frac{\p\tilde H_1}{\p \vfi}\right)}
 =\overline{-\frac{\p }{\p I}\left(\frac{1}{\omega_0}\tilde H_1\frac{\p\tilde H_1}{\p \vfi}\right)}\\
 =&\overline{-\frac{\p }{\p I}\left(\frac{1}{2\omega_0}\frac{\p\tilde H_1^2}{\p \vfi}\right)}=0.
 \end{aligned}
\end{equation}
Combining all the estimates we get the result of the Lemma.

\medskip
{\bf Appendix.  Probability distribution of a pseudophase}

\medskip
Consider some ball  $U$ in the set of initial conditions $D_0\times \Sm^1$. Fix any $r>1$.
According to Proposition \ref{prop1}, solutions of system (\ref{e_system}) with initial conditions in $U\setminus {\cal V}_r$ arrive to the resonance,  and ${ {\rm mes}}\  {\cal V}_r=O(\eps^r)$. For these solutions on can determine phase at the resonance $\vfi_e$ and the corresponding pseudophase $\Xi$ (\ref{eq_Xi}). As phase is determined ${\rm mod}\ 2\pi$, it is natural to consider the pseudophase  ${\rm mod}\ 1$. Thus, consider  the fractional part $\hat \Xi$ of the variable $\Xi$.

Formula (\ref{Xi_a}) shows that $\vfi _e \ {\rm mod } \   2\pi$ and $\hat \Xi$ are very sensitive to change of initial conditions.  A change of $J_0$ of order $\eps$  produces changes  of order 1 in
$\hat\vfi=~\vfi _e \ {\rm mod } \ 2\pi$ and $\hat \Xi$. Therefore, for small $\eps$, it is reasonable to consider $\hat\vfi$ and  $\hat \Xi$ as random values and define their probability distributions.

Denote $U^{\eps}_{({\cal \alpha}, {\cal \beta} )}$ the set of initial points in $U$ such that
$\hat \Xi \in ({\cal \alpha}, {\cal \beta} )\subseteq (0,1)$.

\begin{dfn}\cite{Nei1997}
The value
\begin{equation}
Pr\, (\hat \Xi\in ({\cal \alpha}, {\cal \beta} ))=\lim_{\eps\to 0}\frac{\ {\rm mes} \ U^{\eps}_{({\cal \alpha}, {\cal \beta} )}}{{\rm mes}\ U }
\nonumber
\end{equation}
is called the probability of the event $ \hat \Xi\in ({\cal \alpha}, {\cal \beta} )$.
\end{dfn}
Therefore, the probability of the event $ \hat \Xi\in ({\cal \alpha}, {\cal \beta} )$ the limit (as $\eps \to 0$) of the fraction of  phase volume  in $U$  that is occupied  by  initial conditions
for trajectories with $ \hat \Xi\in ({\cal \alpha}, {\cal \beta} )$. This approach to defining a probability follows  to that  by V.I.Arnold in \cite {Arn_63}.
\begin{proposition}\cite{Nei1997}
The value  $\hat\Xi$ has the
uniform distribution on the interval (0,~ 1):
\begin{equation}
Pr \, (\hat \Xi\in ({\cal \alpha}, {\cal \beta} ))
={\cal \beta} - {\cal \alpha}.
\nonumber
\end{equation}
\end{proposition}

This statement  is justified in \cite{Nei1997}  via computing phase fluxes. It can be obtained as a direct corollary of  formula (\ref{Xi_a}). Moreover, this formula shows, that the  uniform distribution of $\hat \Xi$  arises on each curve  $\{\vfi_0= {\rm const}, y= {\rm const}, x= {\rm const}\}$ in $U$.  This is because  $d\,\hat \Xi_a/dJ_0\sim \eps^{-1}$ when $J_0$ changes along such a curve. Here $\hat \Xi_a$ is the principal term of  the value $\hat \Xi$ given by  (\ref{Xi_a}).  We omit details of the proof.

The probability distribution of $\hat \vfi_e$ can be determined similarly to that of $\hat \Xi$   and obtained from the uniform distribution of $\hat \Xi$. The probability distribution density of  $\hat \vfi_e$ is
\begin{equation}
p(\vfi_e; I_0, y_*,x_*)=   \frac{1}{2\pi}\left( 1+ \frac{1}{b(y_*,x_*)}\frac{\partial H_1(I_0, \vfi_e, y_*,x_*,0)}{\partial \vfi_e}     \right).
\nonumber
\end{equation}
\medskip

Another natural way to define probability distribution  in  problems with a small parameter is suggested by D.V.Anosov (see discussion in \cite {Nei_nonlinearity}).
Fix an initial point  $(I_0, \vfi_0, y_0, x_0)\in U$ and calculate the pseudophase $\hat \Xi$ for this initial point. Note that $\hat \Xi$ depends on $\eps$.  Introduce the set
\begin{equation}
V^{\eps_0}_{({\cal \alpha}, {\cal \beta} )}=\{\eps \in (0, \eps_0): \hat \Xi\in   ({\cal \alpha}, {\cal \beta} )  \}.
\nonumber
\end{equation}

\begin{dfn}
The value
\begin{equation}
Pr'\, (\hat \Xi\in ({\cal \alpha}, {\cal \beta} ))=\lim_{\eps_0\to 0}
\frac{\ {\rm mes} \ V^{\eps_0}_{({\cal \alpha}, {\cal \beta} )}}{\eps_0 }
\nonumber
\end{equation}
is called the probability of the event $ \hat \Xi\in ({\cal \alpha}, {\cal \beta} )$.
\end{dfn}
\begin{proposition}
\label{for_anosov}
Probability $Pr'$ of the event $ \hat \Xi\in ({\cal \alpha}, {\cal \beta} )$ does not depend on the choice of the initial point   $(I_0, \vfi_0, y_0, x_0)\in U$  and is given by the formula
\begin{equation}
Pr' \, (\hat \Xi\in ({\cal \alpha}, {\cal \beta} ))
={\cal \beta} - {\cal \alpha}.
\nonumber
\end{equation}
\end{proposition}
Thus, two definitions of probability lead to the same result. The proof of Proposition \ref{for_anosov} is again based on formula  (\ref{Xi_a}). We omit the proof.
There is an analogous proposition with the complete proof in \cite  {Nei_Okunev_nonlinearity} (see Proposition 3.8 there).

\medskip
{\bf Acknowledgments.} The work was supported by the Leverhulme Trust (Grant No. RPG-2018-143). The authors are thankful to A.V. Artemyev,  S.S. Minkov,  I.S. Shilin and A.A. Vasiliev for
useful discussions.

\vskip 10mm

\noindent Yuyang Gao,

\noindent {\small Department of Mathematical Sciences,}

\noindent {\small Loughborough University, Loughborough LE11 3TU, United Kingdom;}

\noindent {\footnotesize{E-mail : Y.Gao4@lboro.ac.uk}}

\vskip 5mm

\noindent Anatoly Neishtadt,

\noindent {\small Department of Mathematical Sciences,}

\noindent {\small Loughborough University, Loughborough LE11 3TU, United Kingdom;}

\noindent {\footnotesize{E-mail : a.neishtadt@lboro.ac.uk}}

\vskip 5mm

\noindent Alexey Okunev,

\noindent {\small Department of Mathematics,}

\noindent {\small Pennsylvania State University, State College,
Pennsylvania 16802, United States}

\noindent {\footnotesize{E-mail : abo5297@psu.edu}}

\end{document}